\documentclass[final]{siamltex}

\usepackage{epsfig}
\usepackage{amssymb}
\usepackage{amsmath}
\usepackage{amsfonts}
\usepackage{scalefnt}
\usepackage{color}
\usepackage{graphicx}
\usepackage{hyperref}
\usepackage{bookmark}
\usepackage{floatrow}
\newtheorem{remark}{Remark}[section]
\newtheorem{example}{Example}[section]

\title{Joint Spectral Radius and\\ Path-Complete Graph Lyapunov Functions\thanks{A shorter and preliminary version of this work has been
presented in~\cite{HSCC_JSR_Path_complete} as a conference paper.}}

\author{Amir Ali Ahmadi, Rapha\"{e}l Jungers,\\  Pablo A. Parrilo, \and Mardavij
Roozbehani\thanks{Amir Ali Ahmadi, Pablo A. Parrilo, and Mardavij
Roozbehani are with the Laboratory for Information and Decision
Systems, Massachusetts Institute of Technology. Email:
\{\texttt{a\_a\_a}, \texttt{parrilo},
\texttt{mardavij}\}\texttt{@mit.edu.}} Rapha\"{e}l Jungers is an
F.R.S.-FNRS fellow and with the Department of Mathematical
Engineering, Universit\'{e} catholique de Louvain. Email:
\texttt{raphael.jungers@uclouvain.be.} }

\newcommand{\rmjj}[1]{{\color{black} #1}}
\newcommand{\rmj}[1]{{#1}}

\newcommand{\mr}[1]{{\color{black} #1}}

\long\def\aaa#1{{#1}}

\definecolor{Pink}{RGB}{225,0,127}
\newcommand{\amal}[1]{{\color{black} #1}}

\long\def\amral#1{{#1}}

\begin{document}

\maketitle

\begin{abstract}
\noindent We introduce the framework of path-complete graph
Lyapunov functions for approximation of the joint spectral radius.
The approach is based on the analysis of the underlying switched
system via inequalities imposed among multiple Lyapunov functions
associated to a labeled directed graph. Inspired by concepts in
automata theory and symbolic dynamics, we define a class of graphs
called path-complete graphs, and show that any such graph gives
rise to a method for proving stability of the switched system.
This enables us to derive several asymptotically tight hierarchies
of semidefinite programming relaxations that unify and generalize
many existing techniques such as common quadratic, common sum of
squares, \amal{path-dependent quadratic,} and maximum/minimum-of-quadratics Lyapunov functions. We
compare the quality of approximation obtained by certain classes
of path-complete graphs including a family of dual graphs and all
path-complete graphs with two nodes on an alphabet of two
matrices. 
We derive approximation guarantees for several families
of path-complete graphs, such as the De Bruijn graphs. \amal{This provides worst-case perfomance bounds for path-dependent quadratic Lyapunov functions and} a constructive converse Lyapunov
theorem for maximum/minimum-of-quadratics Lyapunov functions.

\end{abstract}



\pagestyle{myheadings} \thispagestyle{plain} \markboth{A. A.
Ahmadi, R. Jungers, P. A. Parrilo, and M. Roozbehani}{Joint
Spectral Radius and Path-Complete Graph Lyapunov Functions}

\section{Introduction}
Given a finite set of square matrices
$\mathcal{A}\mathrel{\mathop:}=\left\{ A_{1},...,A_{m}\right\} $,
their \emph{joint spectral radius} $\rho(\mathcal{A})$ is defined
as
\begin{equation} \rho\left(\mathcal{A}\right)
=\lim_{k\rightarrow\infty}\max_{\sigma
\in\left\{  1,...,m\right\}  ^{k}}\left\Vert A_{\sigma_{k}}...A_{\sigma_{2}%
}A_{\sigma_{1}}\right\Vert ^{1/k},\label{eq:def.jsr}%
\end{equation}
where the quantity $\rho(\mathcal{A})$ is independent of the norm
used in (\ref{eq:def.jsr}). The joint spectral radius (JSR) is a
natural generalization of the spectral radius of a single square
matrix and it characterizes the maximal growth rate that can be
obtained by taking products, of arbitrary length, of all possible
permutations of $A_{1},...,A_{m}$. This concept was introduced by
Rota and Strang~\cite{RoSt60} in the early 60s and has since been
the subject of extensive research within the engineering and the
mathematics communities alike. Aside from a wealth of fascinating
mathematical questions that arise from the JSR, the notion emerges
in many areas of application such as stability of switched linear
dynamical systems, Leontief input-output model of the economy with
uncertain data, computation of the capacity of codes, continuity
of wavelet functions, convergence of consensus algorithms,
trackability of graphs, and many others. See~\cite{Raphael_Book}
and references therein for a recent survey of the theory and
applications of the JSR.

Motivated by the abundance of applications, there has been much
work on efficient computation of the joint spectral radius; see
e.g.~\amal{\cite{grip}},~\cite{BlNT04},~\cite{BlNes05},~\amal{\cite{LeeD06}},~\cite{Pablo_Jadbabaie_JSR_journal},~\amal{\cite{LeeK09}},~\amal{\cite{GZalgorithm}, \cite{GuglielmiZennaro2}, \cite{GP11}} and references therein. Unfortunately, the negative results in the literature certainly restrict the horizon of possibilities.
In~\cite{BlTi2}, Blondel and Tsitsiklis prove that even when the
set $\mathcal{A}$ consists of only two matrices, the question of
testing whether $\rho(\mathcal{A})\leq1$ is undecidable. They also
show that unless P=NP, one cannot compute an approximation
$\hat{\rho}$ of $\rho$ that satisfies
$|\hat{\rho}-\rho|\leq\epsilon\rho$, in a number of steps
polynomial in the bit size of $\mathcal{A}$ and the bit size
of~$\epsilon$~\cite{BlTi3}. It is easy to show that the spectral
radius of any finite product of length $k$ raised to the power of
$1/k$ gives a lower bound on $\rho$~\cite{Raphael_Book}. However,
for reasons that we explain next, our focus will be on computing
upper bounds for $\rho$.

There is an attractive connection between the joint spectral
radius and the stability properties of an arbitrarily switched
linear system; i.e., dynamical systems of the form
\begin{equation}\label{eq:switched.linear.sys}
x_{k+1}=A_{\sigma\left(  k\right)  }x_{k},
\end{equation}
where $\sigma:\mathbb{Z\rightarrow}\left\{  1,...,m\right\}$ is a
map from the set of integers to the set of indices. It is
well-known that $\rho<1$ if and only if system
(\ref{eq:switched.linear.sys}) is \emph{absolutely asymptotically
stable} (AAS), that is, (globally) asymptotically stable for all
switching sequences. Moreover, it is
known~\cite{switched_system_survey} that absolute asymptotic
stability of (\ref{eq:switched.linear.sys}) is equivalent to
absolute asymptotic stability of the linear difference inclusion
\begin{equation}\label{eq:linear.difference.inclusion}
x_{k+1}\in \mbox{co}{\mathcal{A}}\ x_{k},
\end{equation}
where $\mbox{co}{\mathcal{A}}$ here denotes the convex hull of the
set $\mathcal{A}$. Therefore, any method for obtaining upper
bounds on the joint spectral radius provides sufficient conditions
for stability of systems of type (\ref{eq:switched.linear.sys}) or
(\ref{eq:linear.difference.inclusion}). Conversely, if we can
prove absolute asymptotic stability of
(\ref{eq:switched.linear.sys}) or
(\ref{eq:linear.difference.inclusion}) for the set
$\mathcal{A}_\gamma\mathrel{\mathop:}=\{ \gamma
A_{1},\ldots,\gamma A_{m}\}$ for some positive scalar $\gamma$,
then we get an upper bound of $\frac{1}{\gamma}$ on
$\rho(\mathcal{A})$. (This follows from the scaling property of
the JSR: $\rho(\mathcal{A}_\gamma)=\gamma\rho(\mathcal{A})$.) One
advantage of working with the notion of the joint spectral radius
is that it gives a way of rigorously quantifying the performance
guarantee of different techniques for stability analysis of
systems (\ref{eq:switched.linear.sys}) or
(\ref{eq:linear.difference.inclusion}).

Perhaps the most well-established technique for proving stability
of switched systems is the use of a \emph{common (or simultaneous)
Lyapunov function}. The idea here is that if there is a
continuous, positive, and homogeneous (Lyapunov) function
$V:\mathbb{R}^n\rightarrow\mathbb{R}$ that for some $\gamma>1$
satisfies
\begin{equation}
V(\gamma A_ix)\leq V(x) \quad \forall i=1,\ldots,m,\ \ \forall
x\in\mathbb{R}^n,
\end{equation}
(i.e., $V(x)$ decreases no matter which matrix is applied), then
the system in (\ref{eq:switched.linear.sys}) (or in
(\ref{eq:linear.difference.inclusion})) is AAS. Conversely, it is
known that if the system is AAS, then there exists a \emph{convex}
common Lyapunov function (in fact a norm); see e.g.~\cite[p.
24]{Raphael_Book}. However, this function is not in general
finitely constructable. A popular approach has been to try to
approximate this function by a class of functions that we can
efficiently search for using convex optimization and in particular
semidefinite programming. Semidefinite programs (SDPs) can be
solved with arbitrary accuracy in polynomial time and lead to
efficient computational methods for approximation of the JSR. As
an example, if we take the Lyapunov function to be quadratic
(i.e., $V(x)=x^TPx$), then the search for such a Lyapunov function
can be formulated as the following SDP:
\begin{equation}\label{eq:Lyap.CQ.SDP}
\begin{array}{rll}
P&\succ&0 \\
\gamma^2 A_i^TPA_i&\preceq&P \quad \forall i=1,\ldots,m.
\end{array}
\end{equation}

The quality of approximation of common quadratic Lyapunov
functions is a well-studied topic. In particular, it is
known~\cite{BlNT04} that the estimate
$\hat{\rho}_{\mathcal{V}^{2}}$ obtained by this
method\footnote{The estimate $\hat{\rho}_{\mathcal{V}^{2}}$ is the
reciprocal of the largest $\gamma$ that satisfies
(\ref{eq:Lyap.CQ.SDP}) and can be found by bisection.} satisfies
\begin{equation}\label{eq:CQ.bound}
\frac{1}{\sqrt{n}}\hat{\rho}_{\mathcal{V}^{2}}(\mathcal{A})\leq\rho(\mathcal{A})\leq\hat{\rho}_{\mathcal{V}^{2}}(\mathcal{A}),
\end{equation}
where $n$ is the dimension of the matrices. This bound is a direct
consequence of John's ellipsoid theorem and is tight~\cite{Ando98}. \amal{Morover, it is known that applying the common quadratic method to products of increasing length from the set $\mathcal{A}$ gives an asymptotically exact method for the computation of the JSR~\cite{Ando98},~\cite{BlimanFerrari}.}

In~\cite{Pablo_Jadbabaie_JSR_journal}, the use of sum of squares
(SOS) polynomial Lyapunov functions of degree $2d$ was proposed as
a common Lyapunov function for the switched system in
(\ref{eq:switched.linear.sys}). The search for such a Lyapunov
function can again be formulated as a semidefinite program. This
method does considerably better than a common quadratic Lyapunov
function in practice and its estimate
$\hat{\rho}_{\mathcal{V}^{SOS,2d}}$ satisfies the bound
\begin{equation}\label{eq:SOS.bound}
\frac{1}{\sqrt[2d]{\eta}}\hat{\rho}_{\mathcal{V}^{SOS,2d}}(\mathcal{A})\leq\rho(\mathcal{A})\leq\hat{\rho}_{\mathcal{V}^{SOS,2d}}(\mathcal{A}),
\end{equation}
where $\eta=\min\{m,{n+d-1\choose d}\}$. Furthermore, as the
degree $2d$ goes to infinity, the estimate
$\hat{\rho}_{\mathcal{V}^{SOS,2d}}$ converges to the true value of
$\rho$~\cite{Pablo_Jadbabaie_JSR_journal}.

The semidefinite
programming based methods for approximation of the JSR have been
recently generalized and put in the framework of conic
programming~\cite{protasov-jungers-blondel09}.
\amal{We shall also remark that there are powerful techniques for approximation of the JSR that do not use semidefinite programming, such as approaches based on computation of a polytopic norm \cite{GP11}, \cite{GZalgorithm}, \cite{GuglielmiZennaro2}. Research in the computation of the JSR continues to be an active area and each novel technique has the potential to enhance not only our ability to solve certain instances more efficiently, but also our understanding of the relations between the different approaches. An increasing number of the currently available methods for JSR approximation are being (or have been) implemented in the JSR Toolbox, a MATLAB based software package freely available for download~\cite{jsr_toolbox}. Extensive numerical experiments comparing some of the different approaches have been carried out using this toolbox and recently reported in~\cite{Experiment_JSR}.}

\subsection{Contributions and organization}
It is natural to ask whether one can develop better approximation
schemes for the joint spectral radius by using multiple Lyapunov
functions as opposed to requiring simultaneous contractibility of
a single Lyapunov function with respect to all the matrices. More
concretely, our goal is to understand \amal{in what ways} we can write
inequalities among, say, $k$ different Lyapunov functions
$V_1(x),\ldots,V_k(x)$ that imply absolute asymptotic stability of
(\ref{eq:switched.linear.sys}) and can be checked via semidefinite
programming.

The general idea of using several Lyapunov functions for analysis
of switched systems is a very natural one and has already appeared
in the literature (although to our knowledge not in the context of
the approximation of the JSR); see e.g.
\cite{JohRan_PWQ},~\cite{multiple_lyap_Branicky},~\amal{\cite{daafouzbernussou}, \cite{LeeD06}, \cite{LeeK09}, \cite{LeeD06_disturb}, \cite{LeeK08_output}}, \cite{composite_Lyap}, \cite{composite_Lyap2},
\cite{convex_conjugate_Lyap}. Perhaps one of the earliest
references is the work on ``piecewise quadratic Lyapunov
functions'' in~\cite{JohRan_PWQ}. However, this work is in the
different framework of state dependent switching, where the
dynamics switches depending on which region of the space the
trajectory is traversing (as opposed to arbitrary switching). In
this setting, there is a natural way of using several Lyapunov
functions: assign one Lyapunov function per region and ``glue them
together''. Closer to our setting, there is a body of work in the
literature that gives sufficient conditions for existence of
piecewise Lyapunov functions of the type $\max\{x^TP_1x,\ldots,
x^TP_kx\}$, $\min\{x^TP_1x,\ldots,x^TP_kx\}$, and
$\mbox{conv}\{x^TP_1x,\ldots, x^TP_kx\}$, i.e., the pointwise
maximum, the pointwise minimum, and the convex envelope of a set
of quadratic functions \cite{composite_Lyap},
\cite{composite_Lyap2}, \cite{convex_conjugate_Lyap},
\cite{hu-ma-lin}. These works are mostly concerned with analysis
of linear differential inclusions in continuous time, but they
have obvious discrete time counterparts. The main drawback of
these methods is that in their greatest generality, they involve
solving bilinear matrix inequalities, which are non-convex and in
general NP-hard. One therefore has to turn to heuristics, which
have no performance guarantees and their computation time quickly
becomes prohibitive when the dimension of the system increases.
Moreover, \amal{these} methods solely provide sufficient
conditions for stability with no performance guarantees.

\amal{Another body of work which utilizes multiple Lyapunov functions and is of particular interest for us appears in \cite{LeeD06}, \cite{LeeK09}, \cite{LeeD06_disturb}, \cite{LeeK08_output}. In these papers, several fundamental control problems (e.g. stability, feedback stabilizability, detectability, disturbance attenuation, output regulation, etc.) are addressed for discrete-time switched systems using multiple Lyapunov functions and hierarchies of linear matrix inequality (LMI) conditions. The special case of these results that handles the stability question for arbitrarily switched linear systems is directly relevant for our purposes. This includes some of the LMIs associated with the so-called \emph{path-dependent quadratic Lyapunov functions}~\cite{LeeD06}, and another family of LMIs that are in a certain sense dual to those of path-dependent quadratic Lyapunov functions; see~\cite{LeeK09}. In contrast to the piecewise Lyapunov functions discussed previously, these techniques, being naturally SDP-based, do not suffer from computational difficulties associated with solving bilinear matrix inequalities. Moreover, just like the case of sum of squares Lyapunov functions, the hierarchies of LMIs in~\cite{LeeD06}, \cite{LeeK09} are asymptotically exact for computation of the JSR. In other words, the infinite family of the LMIs provides \emph{necessary} and sufficient conditions for switched stability. We will revisit some of these LMIs in this paper, prove approximation guarantees for them, and relate them to common min/max-of-quadratics Lyapunov functions.}


\amal{Motivated by the premise that techniques combining multiple Lyapunov functions and convex optimization provide powerful tools for stability analysis of switched systems, we believe it is important to establish a systematic framework for deriving convex inequalities among multiple Lyapunov functions that imply stability. Moreover, it is naturally desired to understand the performance of the resulting convex programs in terms of approximation of the JSR, just like we do for several classes of common Lyapunov functions (e.g. common quadratic or common SOS). In more concrete terms, the questions that motivate our paper are as follows:} (i) With a focus on conditions that are amenable to convex optimization, what are \amal{\emph{all}} the different ways to write a set of inequalities among $k$ Lyapunov functions that imply absolute asymptotic stability of
(\ref{eq:switched.linear.sys})? Can we give a unifying framework that includes all the previously proposed Lyapunov functions \amal{in the literature}? \amal{Are there new sets of inequalities that have not appeared before?} (ii) Among the different sets of inequalities that imply stability, can we identify some that are more powerful than some other? (iii) The available \amal{(finite) convex programs based on multiple} Lyapunov functions solely provide sufficient conditions for stability with no guarantee on their \amal{approximation quality for the JSR}. Can we give converse theorems that guarantee the existence of a feasible solution to our search for a given accuracy {of approximation}?

\aaa{The contributions of this paper to these questions are as
follows.} We propose a unifying framework based on a
representation of Lyapunov inequalities with labeled graphs and by
making some connections with basic concepts in automata theory.
This is done in Section~\ref{sec:graphs.jsr}, where we define the
notion of a path-complete graph
(Definition~\ref{def:path-complete}) and prove that any such graph
provides an approximation scheme for the JSR
(Theorem~\ref{thm:path.complete.implies.stability}). In
Section~\ref{sec:duality.and.some.families.of.path.complete}, we
give examples of families of path-complete graphs and show that the previously proposed techniques come from particular
classes of path-complete graphs \amal{whose path-completeness is ``easy to detect''} (e.g.,
Corollary~\ref{cor:min.of.quadratics},
Corollary~\ref{cor:max.of.quadratics}, \aaa{and
Remark~\ref{rmk:Lee-Dellerud_Daafouz}}).\footnote{\amal{Although there may be other LMIs in the literature that we are unaware of, it is safe for us to assume that they too must form special cases of our framework. In recent work to be reported elsewhere (see~\cite{Path_complete_converse} for a preliminary version), we have shown that all stability proving Lyapunov inequalities in our setting come from path-complete graphs.}  } \amal{We also show that the concept of path-completeness can easily produce new stability proving LMIs not previously present in the literature (e.g. Proposition~\ref{prop:non.trivial.path.complete.graph} and Remark~\ref{rmk:nontrivial.graph.generalized}).}

In Section~\ref{sec:who.beats.who}, we characterize all the
path-complete graphs with two nodes for the analysis of the JSR of
two matrices.  We \amal{present a full characterization of the partial order induced on these graphs according to their relative performance in approximation of the JSR} (Proposition~\ref{prop:who.beats.who}). In Section~\ref{sec:hscc},
we study in more depth the approximation properties of a
particular pair of ``dual'' path-complete graphs that seem to
perform very well in practice. \amal{The LMIs associated with these dual graphs appear in~\cite{daafouzbernussou}, \cite{LeeD06}, \cite{LeeK09}.} Subsection~\ref{subsec:duality.and.transposition} contains more
general results about duality within path-complete graphs and its
connection to transposition of matrices
(Theorem~\ref{thm:transpose.bound.dual.bound}).
Subsection~\ref{subsec:HSCC.bound} gives an approximation
guarantee for the graphs studied in Section~\ref{sec:hscc}
(Theorem~\ref{thm:HSCC.bound}). \amral{Subsection~\ref{subsec:numerical.examples} contains several numerical
examples, in particular some that come from three application domains: (i) asymptotics of overlap-free words, (ii) computation of the Euler ternary partition function, and (iii) continuity of wavelet functions.} In Section~\ref{sec:converse.thms}, we prove a converse
theorem for the method of max-of-quadratics Lyapunov functions
(Theorem~\ref{thm:converse.max.of.quadratics})  \amal{which tell us how many quadratic Lyapunov functions suffice in worst case to achieve a given approximation quality on the JSR.} \amal{We also derive approximation guarantees for a new class of stability proving LMIs that involve matrix products from the set $\mathcal{A}$ with different lengths} (Theorem~\ref{thm-bound-codes}). Finally, our conclusions and some future directions are presented
in Section~\ref{sec:conclusions.future.directions}.


\section{Path-complete graphs and the joint spectral
radius}\label{sec:graphs.jsr}

In what follows, we will think of the set of matrices
$\mathcal{A}\mathrel{\mathop:}=\left\{ A_{1},...,A_{m}\right\} $
as a finite alphabet and we will often refer to a finite product
of matrices from this set as a \emph{word}. We denote the set of
all words ${A_i}_t\ldots{A_i}_1$ of length $t$ by $\mathcal{A}^t$.
Contrary to the standard convention in automata theory, our
convention is to read a word from right to left. This is in
accordance with the order of matrix multiplication. The set of all
finite words is denoted by $\mathcal{A}^*$; i.e.,
$\mathcal{A}^*=\bigcup\limits_{t\in\mathbb{Z}^+} \mathcal{A}^t$.

\begin{figure}[h]
\centering \scalebox{0.25} {\includegraphics{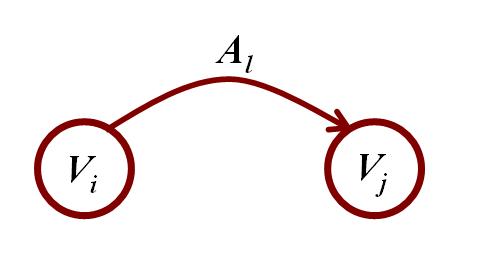}}  
\caption{Graphical representation of Lyapunov inequalities. The
edge in the graph above corresponds to the Lyapunov inequality
$V_j(A_lx)\leq V_i(x)$. Here, $A_l$ can be a single matrix from
$\mathcal{A}$ or a finite product of matrices from $\mathcal{A}$.}
\label{fig:node.arc}
\end{figure}

The basic idea behind our framework is to represent through a
graph all the possible occurrences of products that can appear in
a run of the dynamical system in (\ref{eq:switched.linear.sys}),
and assert via some Lyapunov inequalities that no matter what
occurrence appears, the product must remain stable. A convenient
way of representing these Lyapunov inequalities is via a directed
labeled graph $G(N, E)$. Each node of this graph is associated
with a (continuous, positive definite, and homogeneous) Lyapunov
function $V_i:\mathbb{R}^n\rightarrow\mathbb{R}$, and each edge is
labeled by a finite product of matrices, i.e., by a word from the
set $\mathcal{A}^*$. As illustrated in Figure~\ref{fig:node.arc},
given two nodes with Lyapunov functions $V_i(x)$ and $V_j(x)$ and
an edge going from node $i$ to node $j$ labeled with the matrix
$A_l$, we write the Lyapunov inequality:
\begin{equation}\label{eq:lyap.inequality.rule}
V_j(A_lx)\leq V_i(x) \quad \forall x\in\mathbb{R}^n.
\end{equation}

The problem that we are interested in is to understand which sets
of Lyapunov inequalities imply stability of the switched system in
(\ref{eq:switched.linear.sys}). We will answer this question based
on the corresponding graph.

For reasons that will become clear shortly, we would like to
reduce graphs whose edges have arbitrary labels from the set
$\mathcal{A}^*$ to graphs whose edges have labels from the set
$\mathcal{A}$, i.e., labels of length one. This is explained next.

\begin{definition}\label{def:expanded.graph}
Given a labeled directed graph $G(N,E)$, we define its
\emph{expanded graph} $G^e(N^e,E^e)$ as the outcome of the
following procedure. For every edge $(i,j)\in E$ with label
${A_i}_k\ldots{A_i}_1\in\mathcal{A}^k$, where $k>1$, we remove the
edge $(i,j)$ and replace it with $k$ new edges $(s_q,s_{q+1})\in
E^e\setminus E:\ q\in\{0,\ldots,k-1\}$, where $s_0=i$ and
$s_k=j$.\footnote{It is understood that the node index $s_q$
depends on the original nodes $i$ and $j$. To keep the notation
simple we write $s_q$ instead of $s_{q}^{ij}$.} (These new edges
go from node $i$ through $k-1$ newly added nodes
$s_1,\ldots,s_{k-1}$ and then to node $j$.) We then label the new
edges $(i,s_1),\ldots,(s_q,s_{q+1}),\ldots,(s_{k-1},j)$ with
${A_i}_1,\ldots,{A_i}_k$ respectively.
\end{definition}
\begin{figure}[h]
\centering \scalebox{0.3} {\includegraphics{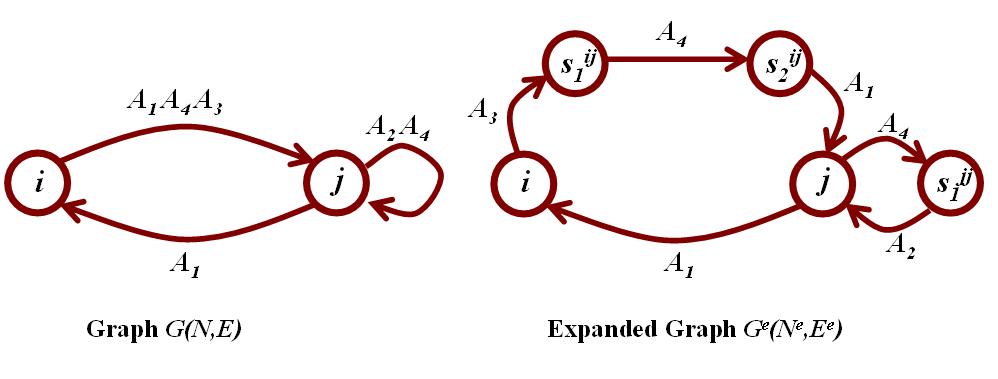}}
\caption{Graph expansion: edges with labels of length more than
one are broken into new edges with labels of length one.}
\label{fig:graph.expansion}
\end{figure}

An example of a graph and its expansion is given in
Figure~\ref{fig:graph.expansion}. Note that if a graph has only
labels of length one, then its expanded graph equals itself. The
next definition is central to our development.
\begin{definition}\label{def:path-complete}
Given a directed graph $G(N,E)$ whose edges are labeled with words
from the set $\mathcal{A}^*$, we say that the graph is
\emph{path-complete}, if for all finite words $A_{\sigma_k}\ldots
A_{\sigma_1}$ of any length $k$ (i.e., for all words in
$\mathcal{A}^*$), there is a directed path in its expanded graph
$G^e(N^e,E^e)$ such that the labels \aaa{on} the edges \aaa{of}
this path are \rmj{the} labels $A_{\sigma_1}$ up to
$A_{\sigma_k}$.
\end{definition}

\begin{figure}[h]
\centering \scalebox{0.4} {\includegraphics{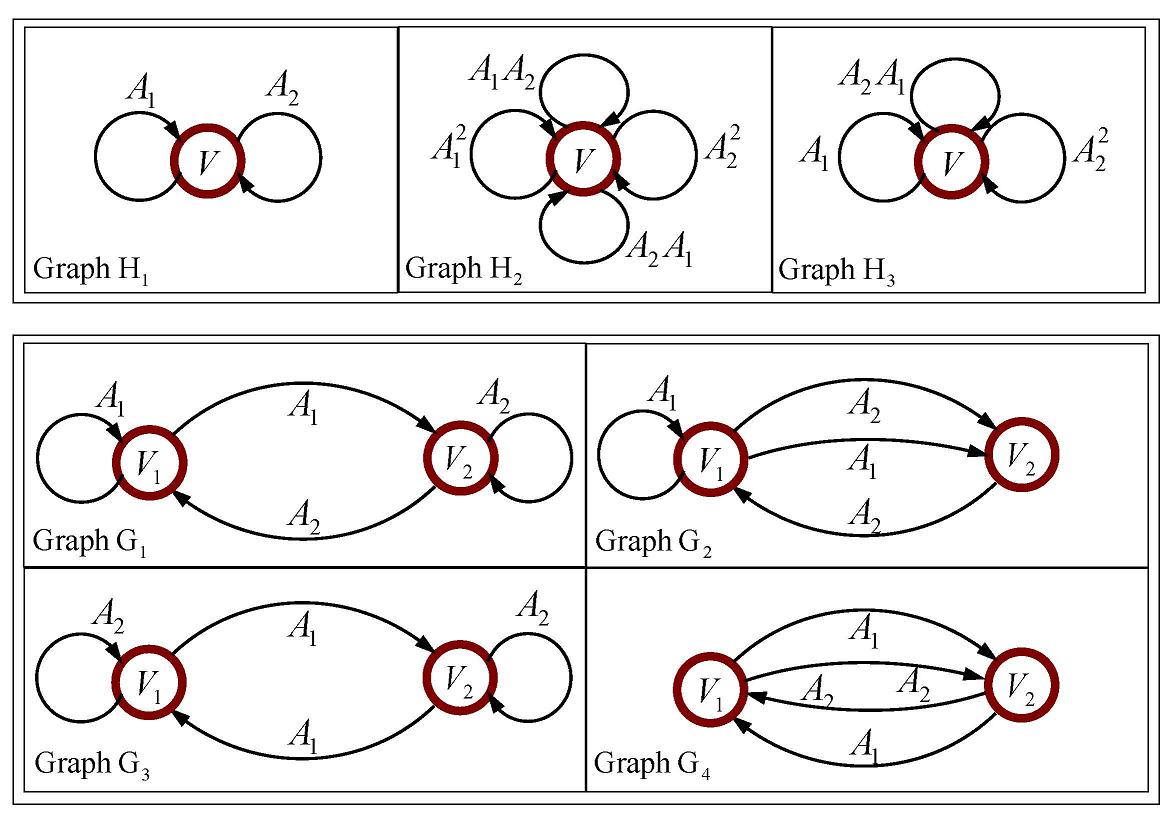}}
\caption{Examples of path-complete graphs for the alphabet
$\{A_1,A_2\}$. If Lyapunov functions satisfying the inequalities
associated with any of these graphs are found, then we get an
upper bound of unity on $\rho(A_1,A_2)$.} \label{fig:jsr.graphs}
\end{figure}

In Figure~\ref{fig:jsr.graphs}, we present seven path-complete
graphs on the alphabet $\mathcal{A}=\{A_1,A_2\}$. The fact that
these graphs are path-complete is easy to see for graphs $H_1,
H_2, G_3,$ and $G_4$, but perhaps not so obvious for graphs $H_3,
G_1,$ and $G_2$. One way to check if a graph is path-complete is
to think of it as a finite automaton by introducing an auxiliary
start node (state) with free transitions to every node and by
making all the other nodes be accepting states. Then, there are
well-known algorithms (see e.g.~\cite[Chap.
4]{Hopcroft_Motwani_Ullman_automata_Book}) that check whether the
language accepted by an automaton is $\mathcal{A}^*$, which is
equivalent to the graph being path-complete. Similar algorithms
exist in the symbolic dynamics literature; see e.g.~\cite[Chap.
3]{Lind_Marcus_symbolic_Book}. Our interest in path-complete
graphs stems from
Theorem~\ref{thm:path.complete.implies.stability} below that
establishes that any such graph gives a method for approximation
of the JSR. We introduce one last definition before we state this
theorem.
%
%

\begin{definition}
Let $\mathcal{A}=\left\{ A_{1},\ldots,A_{m}\right\}  $ be a set of
matrices. Given a path-complete graph $G\left( N,E\right) $ and
$|N|$ functions $V_i(x)$, we say that
$\{V_{i}(x)~|~i=1,\ldots,\left\vert N\right\vert \}$ is a
\emph{graph Lyapunov function (GLF) associated with $G\left(
N,E\right)  $} if
\[
V_{j}\left(  L\left( (i,j)  \right)  x\right)  \leq V_{i}\left(
x\right)  \text{\qquad}\forall x\in \mathbb{R}^n,\ \ \forall\
(i,j)\in E,
\]
where $L\left(  (i,j)  \right) \in\mathcal{A}^{\ast}$ is the label
associated with edge $(i,j)  \in E$ going from node $i$ to node
$j$.
\end{definition}

\begin{theorem}\label{thm:path.complete.implies.stability}
Consider a finite set of matrices
$\mathcal{A}=\{A_1,\ldots,A_m\}$. For a scalar $\gamma>0$, let
$\mathcal{A}_\gamma\mathrel{\mathop:}=\{ \gamma
A_{1},\ldots,\gamma A_{m}\}$. Let $G(N,E)$ be a path-complete
graph whose edges are labeled with words from
$\mathcal{A}_\gamma^*$. If there exist positive, continuous, and
homogeneous\footnote{The requirement of homogeneity can be
replaced by radial unboundedness which is implied by homogeneity
and positivity. However, since the dynamical system in
(\ref{eq:switched.linear.sys}) is homogeneous, there is no
conservatism in asking $V_i(x)$ to be homogeneous.} functions
$V_i(x)$, one per node of the graph, such that
$\{V_{i}(x)~|~i=1,\ldots,\left\vert N\right\vert \}$ is a graph
Lyapunov function associated with $G(N,E)$, then
$\rho(\mathcal{A})\leq\frac{1}{\gamma}$.
\end{theorem}

\begin{proof}
We will first prove the claim for the special case where the edge
labels of $G(N,E)$ belong to $\mathcal{A}_\gamma$ and therefore
$G(N,E)=G^e(N^e,E^e)$. The general case will be reduced to this
case afterwards. Let $d$ be the degree of homogeneity of the
Lyapunov functions $V_i(x)$, i.e., $V_i(\lambda x)=\lambda^d
V_i(x)$ for all $\lambda\in\mathbb{R}$. (The actual value of $d$
is irrelevant.) By positivity, continuity, and homogeneity of
$V_i(x)$, there exist scalars $\alpha_i$ and $\beta_i$ with
$0<\alpha_i\leq\beta_i$ for $i=1,\ldots,|N|$, such that
\begin{equation}\label{eq:homog.bounds}
\alpha_i||x||^d\leq V_i(x)\leq\beta_i||x||^d,
\end{equation}
for all $x\in\mathbb{R}^n$ and for all $i=1,\ldots,|N|$, where
$||x||$ here denotes the Euclidean norm of $x$. Let
\begin{equation}\label{eq:xi}
\xi=\max_{i,j\in\{1,\ldots,|N|\}^2} \frac{\beta_i}{\alpha_j}.
\end{equation}
Now consider an arbitrary product $A_{\sigma_{k}}\ldots
A_{\sigma_{1}}$ of length $k$. Because the graph is path-complete,
there will be a directed path corresponding to this product that
consists of $k$ edges, and goes from some node $i$ to some node
$j$. If we write the chain of $k$ Lyapunov inequalities associated
with these edges (cf. Figure~\ref{fig:node.arc}), then we get
\begin{equation}\nonumber
V_j(\gamma^k A_{\sigma_{k}}\ldots A_{\sigma_{1}}x)\leq V_i(x),
\end{equation}
which by homogeneity of the Lyapunov functions can be rearranged
to
\begin{equation}\label{eq:Vj/Vi.bound}
\left(\frac{V_j(A_{\sigma_{k}}\ldots
A_{\sigma_{1}}x)}{V_i(x)}\right)^{\frac{1}{d}}\leq\frac{1}{\gamma^k}.
\end{equation}
We can now bound the spectral norm of $A_{\sigma_{k}}\ldots
A_{\sigma_{1}}$ as follows:
\begin{eqnarray}\nonumber
||A_{\sigma_{k}}\ldots A_{\sigma_{1}}||&\leq&\max_{x}
\frac{||A_{\sigma_{k}}\ldots A_{\sigma_{1}}x||}{||x||} \\
\nonumber\
&\leq&\left(\frac{\beta_i}{\alpha_j}\right)^{\frac{1}{d}}\max_{x}\frac{V_j^{\frac{1}{d}}(A_{\sigma_{k}}\ldots
A_{\sigma_{1}}x)}{V_i^{\frac{1}{d}}(x)} \\\nonumber \
&\leq&\left(\frac{\beta_i}{\alpha_j}\right)^\frac{1}{d}\frac{1}{\gamma^k}
\\\nonumber
\ &\leq&\xi^{\frac{1}{d}}\frac{1}{\gamma^k},
\end{eqnarray}
where the last three inequalities follow from
(\ref{eq:homog.bounds}), (\ref{eq:Vj/Vi.bound}), and (\ref{eq:xi})
respectively. From the definition of the JSR in
(\ref{eq:def.jsr}), after taking the $k$-th root and the limit
$k\rightarrow\infty$, we get that
$\rho(\mathcal{A})\leq\frac{1}{\gamma}$ and the claim is
established.

Now consider the case where at least one edge of $G(N,E)$ has a
label of length more than one and hence $G^e(N^e,E^e)\neq G(N,E).$\footnote{\amal{A reviewer kindly pointed out an alternative and shorter way of proving
the second part of this theorem, without relying on the notion of expanded graphs.  We
present the proof with expanded graphs because the explicit relationship between the Lyapunov functions of a graph and its expanded version prove to be useful in showing equivalence of certain path-complete graphs in terms of the quality of approximation that they provide on the JSR.}}
We will start with the Lyapunov functions $V_i(x)$ assigned to the
nodes of $G(N,E)$ and from them we will explicitly construct
$|N^e|$ Lyapunov functions for the nodes of $G^e(N^e,E^e)$ that
satisfy the Lyapunov inequalities associated to the edges in
$E^e$. Once this is done, in view of our preceding argument and
the fact that the edges of $G^e(N^e,E^e)$ have labels of length
one by definition, the proof will be completed.


For $j\in N^e$, let us denote the new Lyapunov functions by
$V_j^e(x)$. We give the construction for the case where
$\left\vert N^{e}\right\vert =\left\vert N\right\vert +1.$ The
result for the general case follows by iterating this simple
construction. Let $s\in N^{e}\backslash N$ be the added node in
the expanded graph, and $q,r\in N$ be such that $\left( s,q\right)
\in E^{e}$ and $\left( r,s\right)  \in E^{e}$ with $A_{sq}$ and
$A_{rs}$ as the corresponding labels respectively. Define
\begin{equation}
V_{j}^{e}\left(  x\right)  =\left\{
\begin{array}
[c]{lll}%
V_{j}\left(  x\right)  ,\text{ } & \text{if} & j\in N\medskip\\
V_{q}\left(  A_{sq}x\right)  ,\text{ } & \text{if} & j=s.
\end{array}
\right.  \label{Eqtwo}%
\end{equation}
By construction, $r$ and $q,$ and subsequently, $A_{sq}$ and
$A_{rs}$ are uniquely defined and hence, $\left\{  V_{j}^e\left(
x\right)  ~|~j\in N^{e}\right\}  $ is well defined. We only need
to show that
\begin{align}
V_{q}\left(  A_{sq}x\right)   &  \leq V_{s}^{e}\left(  x\right)
\label{thefirst}\\
V_{s}^{e}\left(  A_{rs}x\right)   &  \leq V_{r}\left(  x\right).
\label{thesecond}%
\end{align}
Inequality (\ref{thefirst}) follows trivially from (\ref{Eqtwo}).
Furthermore,
it follows from (\ref{Eqtwo}) that%
\begin{align*}
V_{s}^{e}\left(  A_{rs}x\right)   &  =V_{q}\left(  A_{sq}A_{rs}x\right)  \\
&  \leq V_{r}\left(  x\right),
\end{align*}
where the inequality follows from the fact that for $i\in N$, the
functions $V_i(x)$ satisfy the Lyapunov inequalities of the edges
of $G\left( N,E\right).$
\end{proof}

\begin{remark}\label{rmk:invertibility}
If the matrix $A_{sq}$ is not invertible, the extended function
$V_{j}^{e}(x)$ as defined in (\ref{Eqtwo}) will only be positive
semidefinite. However, since our goal is to approximate the JSR,
we will never be concerned with invertibility of the matrices in
$\mathcal{A}$. Indeed, since the JSR is continuous in the entries
of the matrices~\cite[p. 18]{Raphael_Book}, we can always perturb
the matrices slightly to make them invertible without changing the
JSR by much. In particular, for any $\alpha>0,$ there exist
$0<\varepsilon, \delta <\alpha$ such that
\[
\hat{A}_{sq}=\frac{A_{sq}+\delta I}{1+\varepsilon}%
\]
is invertible and (\ref{Eqtwo})$-$(\ref{thesecond}) are satisfied
with $A_{sq}=\hat{A}_{sq}.$
\end{remark}

To understand the generality of the framework of ``path-complete
graph Lyapunov funcitons'' more clearly, let us revisit the
path-complete graphs in Figure~\ref{fig:jsr.graphs} for the study
of the case where the set $\mathcal{A}=\{A_1,A_2\}$ consists of
only two matrices. For all of these graphs if our choice for the
Lyapunov functions $V(x)$ or $V_1(x)$ and $V_2(x)$ are quadratic
functions or sum of squares polynomial functions, then we can
formulate the well-established semidefinite programs that search
for these candidate Lyapunov functions.

Graph $H_1$, which is clearly the simplest possible one,
corresponds to the well-known common Lyapunov function approach.
Graph $H_2$ is a common Lyapunov function applied to all products
of length two. This graph also obviously implies
stability.\footnote{By slight abuse of terminology, we say that a
graph implies stability meaning that the associated Lyapunov
inequalities imply stability.} But graph $H_3$ tells us that if we
find a Lyapunov function that decreases whenever $A_1$, $A_2^2$,
and $A_2A_1$ are applied (but with no requirement when $A_1A_2$ is
applied), then we still get stability. This is a priori not
obvious and we believe this approach has not appeared in the
literature before. Graph $H_3$ is also an example that explains \amal{our reasoning behind} the expansion process. Note that for the unexpanded
graph, there is no path for any word of the form $(A_1A_2)^k$ or
of the form $A_2^{2k-1}$, for any $k\in \mathbb{N}.$ However, one
can check that in the expanded graph of graph $H_3$, there is a
path for every finite word, and this in turn allows us to conclude
stability from the Lyapunov inequalities of graph $H_3$.

The remaining graphs in Figure~\ref{fig:jsr.graphs} which all have
two nodes and four edges have a connection to the method of
min-of-quadratics or max-of-quadratics Lyapunov
functions~\cite{composite_Lyap},~\cite{composite_Lyap2},
\cite{convex_conjugate_Lyap}, \cite{hu-ma-lin}. If Lyapunov
inequalities associated with any of these four graphs are
satisfied, then either $\min\{V_1(x),V_2(x)\}$ or
$\max\{V_1(x),V_2(x)\}$ or both serve as a common Lyapunov
function for the switched system. In the next section, we assert
these facts in a more general setting
(Corollaries~\ref{cor:min.of.quadratics}
and~\ref{cor:max.of.quadratics}) and show that these graphs in
some sense belong to ``simplest'' families of path-complete
graphs.

%

\section{Duality and examples of families of path-complete
graphs}\label{sec:duality.and.some.families.of.path.complete}

Now that we have shown that \emph{any} path-complete graph yields
a method for proving stability of switched systems, our next focus
is naturally on showing how one can produce graphs that are
path-complete. Before we proceed to some basic constructions of
such graphs, let us define a notion of \emph{duality} among graphs
which essentially doubles the number of path-complete graphs that
we can generate.

\begin{definition}\label{def:dual.graph}
Given a directed graph $G(N,E)$ whose edges are labeled with words
in $\mathcal{A}^*$, we define its \emph{dual graph} $G'(N,E')$ to
be the graph obtained by reversing the direction of the edges of
$G$, and changing the labels $A_{\sigma_k}\ldots A_{\sigma_1}$ of
every edge of $G$ to its reversed version $A_{\sigma_1}\ldots
A_{\sigma_k}$.
\end{definition}

\begin{figure}[h]
\centering \scalebox{0.09} {\includegraphics{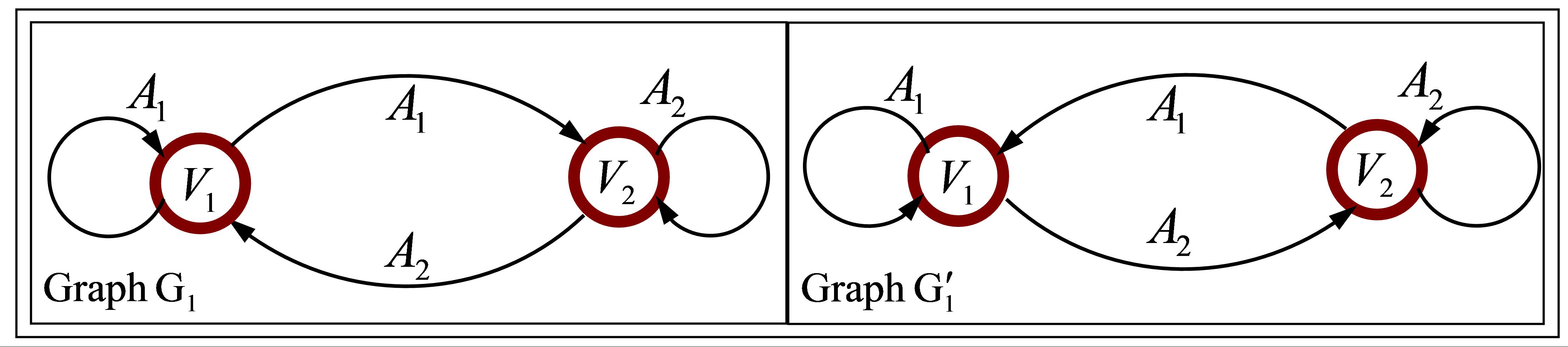}}  
\caption{An example of a pair of dual graphs.}
\label{fig:dual.graphs}
\end{figure}

An example of a pair of dual graphs with labels of length one is
given in Figure~\ref{fig:dual.graphs}. The following theorem
relates dual graphs and path-completeness.

\begin{theorem}\label{thm:path.complete.dual.path.complete}
If a graph $G(N,E)$ is path-complete, then its dual graph
$G'(N,E')$ is also path-complete.
\end{theorem}

\begin{proof}
Consider an arbitrary finite word $A_{i_k}\ldots A_{i_1}$. By
definition of path-completeness, our task is to show that there
exists a path corresponding to this word in the expanded graph of
the dual graph $G'$. It is easy to see that the expanded graph of
the dual graph of $G$ is the same as the dual graph of the
expanded graph of $G$; i.e,
$G'^e(N^e,E'^e)=G^{e^{'}}(N^e,E^{e^{'}})$. Therefore, we show a
path for $A_{i_k}\ldots A_{i_1}$ in $G^{e^{'}}$. Consider the
reversed word $A_{i_i}\ldots A_{i_k}$. Since $G$ is path-complete,
there is a path corresponding to this reversed word in $G^e$. Now
if we just trace this path backwards, we get exactly a path for
the original word $A_{i_k}\ldots A_{i_1}$ in $G^{e^{'}}$. This
completes the proof.
\end{proof}

The next proposition offers a very simple construction for
obtaining a large family of path-complete graphs with labels of
length one.

\begin{proposition}\label{prop:in.out.going.edges.path.complete}
A graph having any of the two properties below is path-complete.

Property (i): every node has outgoing edges with all the labels in
$\mathcal{A}$.

Property (ii): every node has incoming edges with all the labels
in $\mathcal{A}$.
\end{proposition}

\begin{proof}
If a graph has Property (i), then it is obviously path-complete.
If a graph has Property (ii), then its dual has Property (i) and
therefore by Theorem~\ref{thm:path.complete.dual.path.complete} it
is path-complete.
\end{proof}

Examples of path-complete graphs that fall in the category of this
proposition include graphs $G_1,G_2,G_3,$ and $G_4$ in
Figure~\ref{fig:jsr.graphs} and all of their dual graphs. By
combining the previous proposition with
Theorem~\ref{thm:path.complete.implies.stability}, we obtain the
following two simple corollaries which unify several linear matrix
inequalities (LMIs) that have been previously proposed in the
literature. These corollaries also provide a link to
min/max-of-quadratics Lyapunov functions. Different special cases
of these LMIs have appeared
in~\cite{composite_Lyap},~\cite{composite_Lyap2},
\cite{convex_conjugate_Lyap}, \cite{hu-ma-lin}, \cite{LeeD06},
\cite{daafouzbernussou}, \amal{\cite{LeeK09}}. Note that the framework of path-complete
graph Lyapunov functions makes the proof of the fact that these
LMIs imply stability immediate. \amal{We also remark that the following corollaries, and hence the graphs in Proposition~\ref{prop:in.out.going.edges.path.complete}, already include infinite subsets of path-complete graphs that are not only sufficient for stability of (\ref{eq:switched.linear.sys}), but also necessary. Examples of such infinite sets of LMIs with their proofs of necessity are given in \cite{LeeD06}, \amal{\cite{LeeK09}}.}

\begin{corollary}\label{cor:min.of.quadratics}
Consider the set $\mathcal{A}=\{A_1,\ldots,A_m\}$ and the
associated switched linear system in
(\ref{eq:switched.linear.sys}) or
(\ref{eq:linear.difference.inclusion}). If there exist \amal{$K$}
positive definite matrices $P_j$ such that
\begin{eqnarray}\label{eq:min.quadratics.LMIs}
\amal{\forall (i,k)\in\{1,\ldots,m\}\times\{1,\ldots,K\}},\ \exists j\in\{1,\ldots,\amal{K}\}\
\nonumber
\\
\mbox{such that}\quad \quad \gamma^2A_i^TP_jA_i\preceq P_k,
\end{eqnarray}
for some $\gamma>1$, then the system is absolutely asymptotically
stable, i.e., $\rho(\mathcal{A})<1$. Moreover, the pointwise
minimum
$$\min\{x^TP_1x,\ldots,\amal{x^TP_Kx}\}$$ of the quadratic functions serves
as a common Lyapunov function.
\end{corollary}

\begin{proof}
The inequalities in (\ref{eq:min.quadratics.LMIs}) imply that
every node of the associated graph has outgoing edges labeled with
all the different $m$ matrices. Therefore, by
Proposition~\ref{prop:in.out.going.edges.path.complete} the graph
is path-complete, and by
Theorem~\ref{thm:path.complete.implies.stability} this implies
absolute asymptotic stability. The proof that the pointwise
minimum of the quadratics is a common Lyapunov function is easy
and left to the reader.
\end{proof}

\begin{corollary}\label{cor:max.of.quadratics}
Consider the set $\mathcal{A}=\{A_1,\ldots,A_m\}$ and the
associated switched linear system in
(\ref{eq:switched.linear.sys}) or
(\ref{eq:linear.difference.inclusion}). If there exist \amal{$K$}
positive definite matrices $P_j$ such that
\begin{eqnarray}\label{eq:max.quadratics.LMIs}
\amal{\forall (i,j)\in\{1,\ldots,m\}\times\{1,\ldots,K\}},\ \exists k\in\{1,\ldots,\amal{K}\}\
\nonumber
\\
\mbox{such that}\quad \quad \gamma^2A_i^TP_jA_i\preceq P_k,
\end{eqnarray}
for some $\gamma>1$, then the system is absolutely asymptotically
stable, i.e., $\rho(\mathcal{A})<1$. Moreover, the pointwise
maximum
$$\max\{x^TP_1x,\ldots,\amal{x^TP_Kx}\}$$ of the quadratic functions serves
as a common Lyapunov function.
\end{corollary}

\begin{proof}
The inequalities in (\ref{eq:max.quadratics.LMIs}) imply that
every node of the associated graph has incoming edges labeled with
all the different $m$ matrices. Therefore, by
Proposition~\ref{prop:in.out.going.edges.path.complete} the graph
is path-complete and the proof of absolute asymptotic stability
then follows. The proof that the pointwise maximum of the
quadratics is a common Lyapunov function is again left to the
reader.
\end{proof}

\begin{remark} The linear matrix inequalities in (\ref{eq:min.quadratics.LMIs}) and (\ref{eq:max.quadratics.LMIs}) are (convex) sufficient
conditions for existence of min-of-quadratics or max-of-quadratics
Lyapunov functions. The converse is not true. The works in
~\cite{composite_Lyap},~\cite{composite_Lyap2},
\cite{convex_conjugate_Lyap}, \cite{hu-ma-lin} have additional
multipliers in (\ref{eq:min.quadratics.LMIs}) and
(\ref{eq:max.quadratics.LMIs}) that make the inequalities
non-convex but when solved with a heuristic method contain a
larger family of min-of-quadratics and max-of-quadratics Lyapunov
functions. Even if the non-convex inequalities with multipliers
could be solved exactly, except for special cases where the
$\mathcal{S}$-procedure is exact (e.g., the case of two quadratic
functions), these methods still do not completely characterize
min-of-quadratics and max-of-quadratics functions.
%
\end{remark}

\aaa{\begin{remark} \label{rmk:Lee-Dellerud_Daafouz} \amal{The LMIs associated with ``path-dependent quadratic Lyapunov functions'' of any given path length (see~\cite{LeeD06}) and the LMIs associated with ``parameter dependent
Lyapunov functions''~\cite{daafouzbernussou}---when specialized to the analysis of
arbitrarily switched linear systems---are special cases of
Corollary~\ref{cor:min.of.quadratics}
and~\ref{cor:max.of.quadratics} respectively}. This observation
makes a connection between these techniques and
min/max-of-quadratics Lyapunov functions which is not established
in~\cite{LeeD06},~\cite{daafouzbernussou}. It is also interesting
to note that the path-complete graph corresponding to \amal{the LMIs of path-dependent quadratic Lyapunov functions of any path length (see Theorem 9 in~\cite{LeeD06}) is the well-known
De Bruijn graph~\cite{GraphTheory_Handbook}. The ``path length'' of these Lyapunov functions is interestingly the dimension of the De Bruijn graph}. We will analyze the bound on the JSR obtained by analysis via this path-complete graph
in later sections \amal{since we have empirically observed that path-dependent quadratic Lyapunov functions are among the most powerful ones in comparison to all of our graphs}. 

\end{remark}}

\begin{figure}[h]
\centering \scalebox{0.5}
{\includegraphics{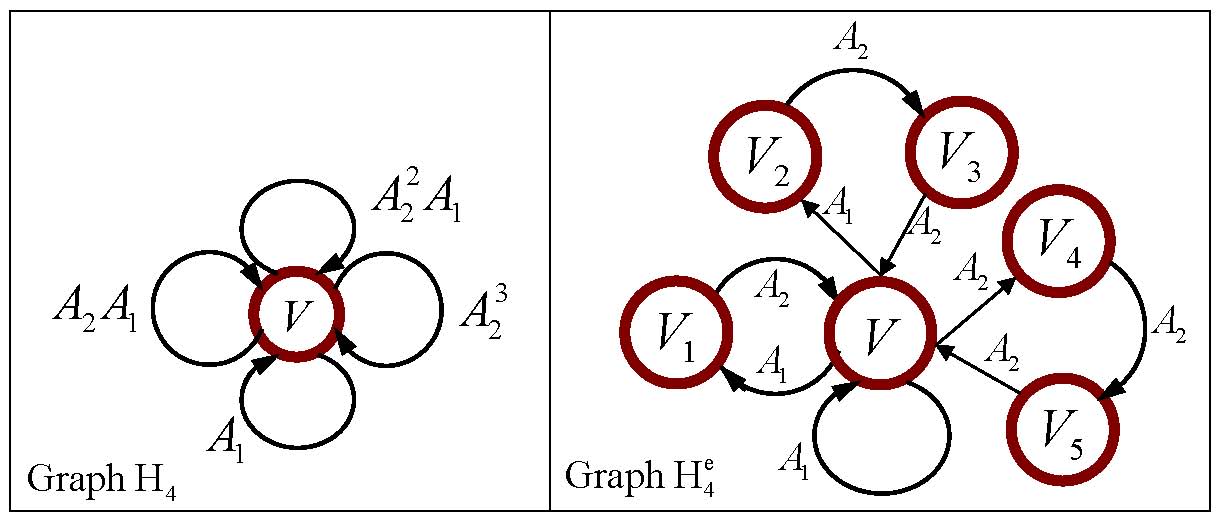}} \caption{The
path-complete graphs corresponding to
Proposition~\ref{prop:non.trivial.path.complete.graph}.}
\label{fig:non.trivial.path.complete}
\end{figure}

The set of path-complete graphs is much broader than \amal{the} family of graphs constructed in
Proposition~\ref{prop:in.out.going.edges.path.complete}. Indeed,
there are many graphs that are path-complete without having
outgoing (or incoming) edges with all the labels on every node;
see e.g. graph $H_4^e$ in
Figure~\ref{fig:non.trivial.path.complete}. This in turn means
that there are several interesting and unexplored Lyapunov
inequalities that we can impose for proving stability of switched
systems. Below, we give one particular example of such
``non-obvious'' inequalities for the case of switching between two
matrices.

\begin{proposition}\label{prop:non.trivial.path.complete.graph}
Consider the set $\mathcal{A}=\{A_1,A_2\}$ and the switched linear
system in (\ref{eq:switched.linear.sys}) or
(\ref{eq:linear.difference.inclusion}). If there exist a positive
definite matrix $P$ such that
\begin{eqnarray}\label{eq:non.trivial.path.complete.graph}
\gamma^2A_1^TPA_1\preceq P, \nonumber \\
\gamma^4(A_2A_1)^TP(A_2A_1)\preceq P, \nonumber \\
\gamma^6(A_2^2A_1)^TP(A_2^2A_1)\preceq P, \nonumber \\
\gamma^6A_2^{3^T}PA_2^3\preceq P, \nonumber
\end{eqnarray}
for some $\gamma>1$, then the system is absolutely asymptotically
stable, i.e., $\rho(\mathcal{A})<1$.
\end{proposition}
\begin{proof}
The graph $H_4$ associated with the LMIs above and its expanded
version $H_4^e$ are drawn in
Figure~\ref{fig:non.trivial.path.complete}. We leave it as an
exercise for the reader to show (e.g. by induction on the length
of the word) that there is path for every finite word in $H_4^e$.
Therefore, $H_4$ is path-complete and in view of
Theorem~\ref{thm:path.complete.implies.stability} the claim is
established.
\end{proof}

\begin{remark}\label{rmk:nontrivial.graph.generalized}
Proposition~\ref{prop:non.trivial.path.complete.graph} can be
generalized as follows: If a single Lyapunov function decreases
with respect to the matrix products
$$\{A_1,A_2A_1,A_2^2A_1,\ldots,A_2^{k-1}A_1,A_2^k\}$$ for some
integer $k\geq 1$, then $\rho(A_1,A_2)<1$. We omit the proof of
this generalization due to space limitations. We will later prove
(Theorem~\ref{thm-bound-codes}) a bound for the quality of
approximation of path-complete graphs of this type, where a common
Lyapunov function is required to decrease with respect to products
of different lengths.
\end{remark}

When we have so many different ways of imposing conditions for
stability, it is natural to ask which ones are more powerful. The
answer clearly depends on the combinatorial structure of the
graphs and does not seem to be easy in general. Nevertheless, in
the next section, we compare the performance of all path-complete
graphs with two nodes for analysis of switched systems with two
matrices. Some interesting connections between the bounds obtained
from these graphs will arise. For example, we will see that the
graphs $H_1, G_3,$ and $G_4$ always give the same bound on the
joint spectral radius; i.e., one graph will succeed in proving
stability if and only if the other two will. So, there is no point
in increasing the number of decision variables and the number of
constraints and impose $G_3$ or $G_4$ in place of $H_1$. The same
is true for the graphs in $H_3$ and $G_2$, which makes graph $H_3$
preferable to graph $G_2$. (See
Proposition~\ref{prop:who.beats.who}.)

\section{Path-complete graphs with two nodes}\label{sec:who.beats.who}

In this section, we characterize the set of all path-complete
graphs consisting of two nodes, an alphabet set
$\mathcal{A}=\{A_{1},A_{2}\},$ and edge labels of unit length. We
will elaborate on the set of all admissible topologies arising in
this setup and compare the performance---in the sense of
conservatism of the ensuing analysis---of different path-complete
graph topologies.

Before we proceed, we introduce a notation that will prove to be
convenient in Subsection~\ref{subsec:comparison.of.who.beats.who}:
Given a labeled graph $G(N,E)$ associated with two matrices $A_1$
and $A_2$, we denote by $\overline{G}(N,E)$, the graph obtained by
swapping of $A_1$ and $A_2$ in all the labels on every edge.

\subsection{The set of path-complete graphs}




The next lemma establishes that for thorough analysis of the case
of two matrices and two nodes, we only need to examine graphs with
four or fewer edges.

\begin{lemma}
\label{lessthan5}Let $G\left(  \left\{  1,2\right\}  ,E\right)  $
be a path-complete graph with labels of length one for
$\mathcal{A}=\{A_{1},A_{2}\}$. Let $\left\{ V_{1},V_{2}\right\}$
be a graph Lyapunov function for $G.$ If $\left\vert E\right\vert
>4,$ then, either
\newline\hspace*{0.15in}(i) there exists $\hat{e}\in E$ such that $G\left(  \left\{  1,2\right\}, E\backslash \hat{e}\right)  $ is a path-complete graph,
\newline\hspace*{0.15in} or
\newline\hspace*{0.15in}(ii) either $V_{1}$ or $V_{2}$ or both are common Lyapunov functions for {$\mathcal{A}.$}
\end{lemma}

\begin{proof}
If $\left\vert E\right\vert >4,$ then at least one node has three
or more outgoing edges. Without loss of generality let node $1$ be
a node with exactly three outgoing edges $e_{1},e_{2},e_{3}$, and
let $L\left( e_{1}\right)  =L\left( e_{2}\right)  =A_{1}.$ Let
$\mathcal{D}\left(  e\right)  $ denote the destination node of an
edge $e\in E.$ If $\mathcal{D}\left(
e_{1}\right)  =\mathcal{D}\left(  e_{2}\right),$ then $e_{1}$ (or $e_{2}%
$)\ can be removed without changing the output set of words. If
$\mathcal{D}\left( e_{1}\right)  \neq\mathcal{D}\left(
e_{2}\right)  ,$ assume, without loss of generality, that
$\mathcal{D}\left(  e_{1}\right)  =1$ and $\mathcal{D}\left(
e_{2}\right)  =2.$ Now, if $L\left(  e_{3}\right)  =A_{1},$ then
regardless of its destination node, $e_{3}$ can be removed. If
$L\left(  e_{3}\right)  =A_{2}$ and $\mathcal{D}\left(
e_{3}\right) =1$, then $V_1$ is a common Lyapunov function for
$\mathcal{A}$. The only remaining possibility is that $L\left(
e_{3}\right)  =A_{2}$ and $\mathcal{D}\left(  e_{3}\right) =2. $
Note that there must be an edge $e_4 \in E$ from node $2$ to node
$1$, otherwise either node $2$ would have two self-edges with the
same label or $V_2$ would be a common Lyapunov function for
$\mathcal{A}$. If $L(e_4)=A_2$ then it can be verified that
$G(\{1,2\},\{e_1,e_2,e_3,e_4\})$ is path-complete and thus all
other edge can be removed. If there is no edge from node $2$ to
node $1$ with label $A_2$ then $L(e_4)=A_1$ and node $2$ must have
a self-edge $e_5 \in E$ with label $L(e_5)=A_2$, otherwise the
graph would not be path-complete. In this case, it can be verified
that $e_{2}$ can be removed without affecting the output set of
words.
\end{proof}

One can easily verify that a path-complete graph with two nodes
and fewer than four edges must necessarily place two self-loops
with different labels on one node, which necessitates existence of
a common Lyapunov function for the underlying switched system.
Since we are interested in exploiting the favorable properties of
\aaa{graph Lyapunov functions} in approximation of the JSR, we
will focus on graphs with four edges.


\subsection{Comparison of performance}\label{subsec:comparison.of.who.beats.who}

It can be verified that for path-complete graphs with two nodes,
four edges, and two matrices, and without multiple self-loops on a
single node, there are a total of nine distinct graph topologies
to consider. Of the nine graphs, six have the property that every
node has two incoming edges with different labels. These are
graphs $G_1,~G_2,~\overline{G}_2,~G_3,~\overline{G}_3,$ and $G_4$
(Figure \ref{fig:jsr.graphs}). Note that $\overline{G}_1=G_1$ and
$\overline{G}_4=G_4$. The duals of these six graphs, i.e.,
$G_1^\prime,~G_2^\prime,~\overline{G}_2^\prime,~G_3^\prime=G_3,~\overline{G}_3^\prime=\overline{G}_3,$
and $G_4^\prime=G_4$ have the property that every node has two
outgoing edges with different labels. Evidently,
$G_3,~\overline{G}_3,$ and $G_4$ are \emph{self-dual graphs},
i.e., they are isomorphic to their dual graphs. The self-dual
graphs are least interesting to us since, as we will show, they
necessitate existence of a common Lyapunov function for
$\mathcal{A}$ (cf. Proposition \ref{prop8}, equation
(\ref{theselfduals}))$.$


Note that all of these graphs perform at least as well as a common
Lyapunov function because we can always take $V_{1}\left( x\right)
=V_{2}\left( x\right)  $. Furthermore, we know from
Corollaries~\ref{cor:max.of.quadratics}
and~\ref{cor:min.of.quadratics} that if Lyapunov inequalities
associated with $G_1,~G_2,~\overline{G}_2,~G_3,~\overline{G}_3,$
and $G_4$ are satisfied, then $\max\left\{ V_{1}\left( x\right)
,V_{2}\left( x\right) \right\} $ is a common Lyapunov function,
whereas, in the case of graphs
$G_1^\prime,~G_2^\prime,~\overline{G}_2^\prime,~G_3^\prime,~\overline{G}_3^\prime$,
and $G_4^\prime$, the function  $\min\left\{ V_{1}\left( x\right)
,V_{2}\left( x\right) \right\} $ would serve as a common Lyapunov
function. Clearly, for the self-dual graphs $G_3,~\overline{G}_3,$
and $G_4$ both $\max\left\{  V_{1}\left( x\right) ,V_{2}\left(
x\right) \right\}$ and $\min\left\{ V_{1}\left(  x\right)
,V_{2}\left( x\right) \right\}$ are common Lyapunov functions.

\textbf{Notation:} Given a set of matrices $\mathcal{A}=\left\{  A_{1}%
,\ldots,A_{m}\right\}  ,$ a path-complete graph $G\left(
N,E\right) ,$ and a
class of functions $\mathcal{V},$ we denote by $\hat{\rho}_{\mathcal{V}}%
,_{G}\left(  \mathcal{A}\right)  ,$ the upper bound on the JSR of
$\mathcal{A}$ that can be obtained by numerical optimization of
GLFs $V_{i}\in \mathcal{V},~i\in N  ,$ defined over $G.$ With a
slight abuse of notation, we denote by $\hat{\rho
}_{\mathcal{V}}\left( \mathcal{A}\right)  ,$ the upper bound that
is obtained by using a common Lyapunov function $V\in\mathcal{V}.$

\aaa{

\begin{proposition}
\label{prop:who.beats.who} \label{prop8} Consider the set
$\mathcal{A}=\left\{A_{1},A_{2}\right\}  ,$  and let
$G_1,~G_2,~G_3,~G_4$, and $H_3$ be the path-complete graphs shown
in Figure \ref{fig:jsr.graphs}. Then, the upper bounds on the JSR
of $\mathcal{A}$ obtained via the associated GLFs satisfy the
following relations:
\begin{equation}
\hat{\rho}_{\mathcal{V}},_{G_{1}}\left(  \mathcal{A}\right)
=\hat{\rho }_{\mathcal{V}},_{G_{1}^{\prime}}\left(
\mathcal{A}\right)
\label{hscceqdual}%
\end{equation}
and
\begin{equation}
\hat{\rho }_{\mathcal{V}}\left(
\mathcal{A}\right)=\hat{\rho}_{\mathcal{V}},_{G_{3}}\left(
\mathcal{A}\right)
=\hat{\rho}_{\mathcal{V}},_{\overline{G}_{3}}\left(
\mathcal{A}\right) =\hat{\rho}_{\mathcal{V}},_{G_{4}}\left(
\mathcal{A}\right)
\label{theselfduals}%
\end{equation}
and
\begin{equation}
\hat{\rho}_{\mathcal{V}},_{G_{2}}\left(  \mathcal{A}\right)
=\hat{\rho }_{\mathcal{V}},_{{{H}_3}}\left(\mathcal{A}\right)
,\text{\qquad}\hat{\rho}_{\mathcal{V}},_{\overline{G}_2}\left(\mathcal{A}\right)
=\hat{\rho}_{\mathcal{V}},_{{\overline{{H}}_3}}\left(\mathcal{A}\right)
\label{theprimals}%
\end{equation}
and
\begin{equation}
\hat{\rho}_{\mathcal{V}},_{G_{2}^{\prime}}\left(
\mathcal{A}\right)
=\hat{\rho}_{\mathcal{V}},_{{{H}_3^{\prime}}}\left(
\mathcal{A}\right)  ,\text{\qquad}\hat{\rho}_{\mathcal{V}},_{\overline{G}_2^{\prime}%
}\left(  \mathcal{A}\right)  =\hat{\rho}_{\mathcal{V}},_{{\overline{{H}}_3^{\prime}}}\left(  \mathcal{A}\right).  \label{theduals}%
\end{equation}
\end{proposition}
}

\begin{proof}
A proof of \eqref{hscceqdual} in more generality is provided in
Section \ref{sec:hscc} (cf. Corollary
\ref{cor:HSCC.invariance.under.transpose}). The proof of
(\ref{theselfduals}) is based on symmetry arguments. Let
$\left\{V_{1},V_{2}\right\}$ be a GLF associated with $G_{3}$
($V_1$ is associated with node $1$ and $V_2$ is associated with
node $2$). Then, by symmetry, $\left\{V_{2},V_{1}\right\}$ is also
a GLF for $G_{3}$ (where $V_1$ is associated with node $2$ and
$V_2$ is associated with node $1$). Therefore, letting
$V=V_{1}+V_{2}$, we have that $\left\{V,V\right\}$ is a GLF for
$G_{3}$ and thus, $V=V_1+V_2$ is also a common Lyapunov function
for $\mathcal{A},$ which implies that
$\hat{\rho}_{\mathcal{V}},_{G_{3}}\left( \mathcal{A}\right)
\geq\hat{\rho }_{\mathcal{V}}\left( \mathcal{A}\right)  .$ The
other direction is trivial: If $V\in\mathcal{V}$ is a common
Lyapunov function for $\mathcal{A},$ then $\left\{
V_{1},V_{2}~|~V_{1}=V_{2}=V\right\}  $ is a GLF associated with
$G_{3},$ and hence, $\hat{\rho}_{\mathcal{V}},_{G_{3}}\left(  \mathcal{A}%
\right)  \leq\hat{\rho}_{\mathcal{V}}\left(  \mathcal{A}\right) .$
Identical arguments based on symmetry hold for
${\overline{G}_{3}}$ and ${G_{4}}$. We now prove the left equality
in (\ref{theprimals}), the proofs for the remaining equalities in
(\ref{theprimals}) and (\ref{theduals}) are analogous. The
equivalence between $G_2$ and $H_3$ is a special case of the
relation between a graph and its \emph{reduced} model, obtained by
removing a node without any self-loops, adding a new edge per each
pair of incoming and outgoing edges to that node, and then
labeling the new edges by taking the composition of the labels of
the corresponding incoming and outgoing edges in the original
graph; see~\cite{Roozbehani2008},~\cite[Chap.
5]{MardavijRoozbehani2008}. Note that $H_3$ is an offspring of
$G_2$ in this sense. This intuition helps construct a proof. Let
$\left\{ V_{1},V_{2}\right\}  $ be a GLF associated with $G_{2}.$
It can be verified that $V_{1}$ is a Lyapunov function associated
with ${{H}_3},$ and therefore,
$\hat{\rho}_{\mathcal{V}},_{{{H}_3}}\left(\mathcal{A}\right)
\leq\hat{\rho}_{\mathcal{V}},_{G_{2}}\left(\mathcal{A}\right) .$
Similarly, if $V\in\mathcal{V}$ is a Lyapunov function associated
with ${{H}_3},$ then one can check that $\left\{
V_{1},V_{2}~|~V_{1}\left(  x\right)  = V\left( x\right)
,V_{2}\left(  x\right)  =V\left(  A_{2}x\right) \right\} $ is a
GLF associated with $G_{2},$ and hence,
$\hat{\rho}_{\mathcal{V}},_{{{H}_3 }}\left(  \mathcal{A}\right)
\geq\hat{\rho}_{\mathcal{V}},_{G_{2}}\left( \mathcal{A}\right).$
\end{proof}


\begin{figure}[h]
\centering \scalebox{0.25} {\includegraphics{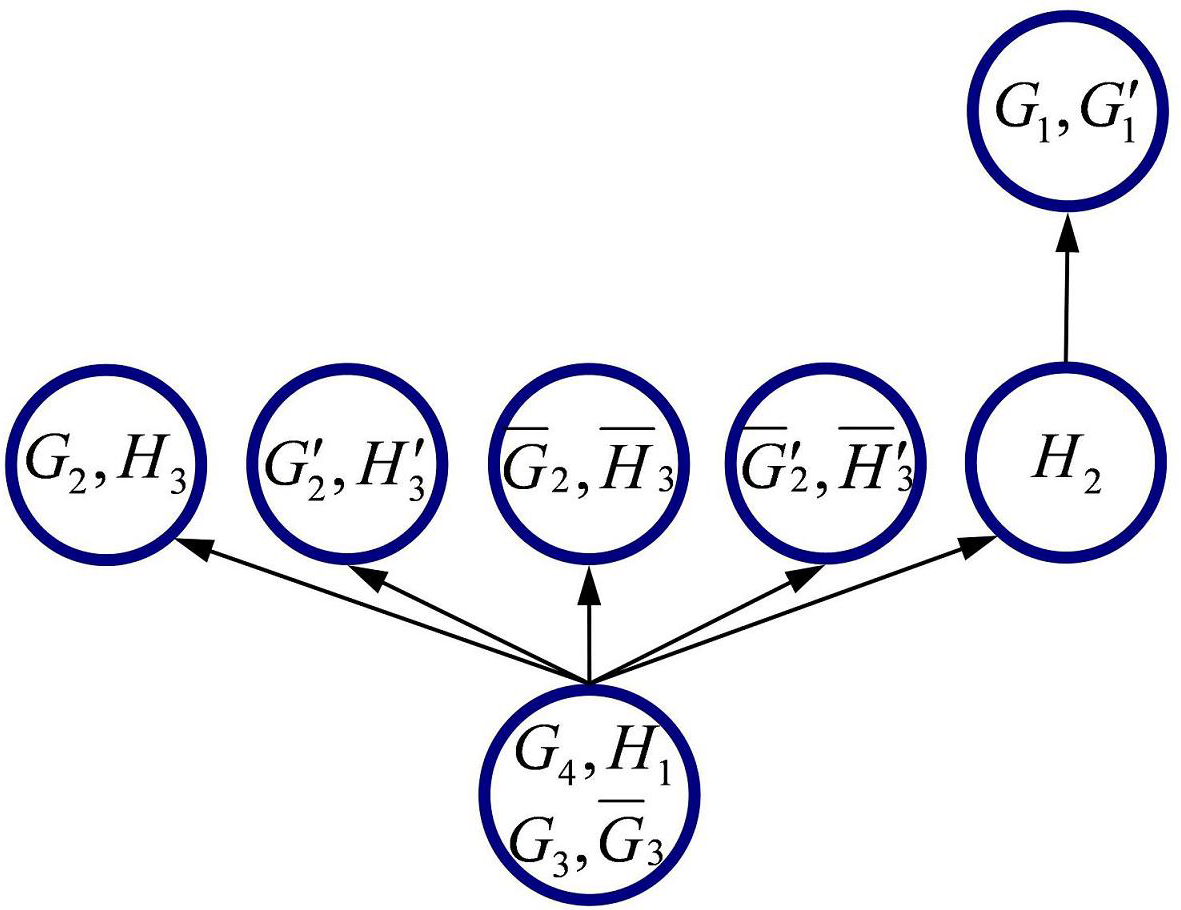}}
\caption{\aaa{A \amal{Hasse} diagram describing the relative performance of the
path-complete graphs of Figure~\ref{fig:jsr.graphs} together with
their duals and label permutations. The graphs placed in the same
circle always give the same approximation of the JSR. A graph at
the end of an arrow results in an approximation of the JSR that is
always at least as good as that of the graph at the start of the
arrow. When there is no directed path between two graphs in this
diagram, either graph can outperform the other depending on the
set of matrices $\mathcal{A}$.}} \label{fig:Hasse.diag}
\end{figure}

\begin{remark}
\label{unverifiedremark}Proposition \ref{prop8} (equation
\ref{hscceqdual}) establishes the equivalence of the bounds
obtained from the pair of dual graphs $G_{1}$ and
$G_{1}^{\prime}$. This, however, is not true for graphs $G_{2}$
and $\overline{G}_2$ as there exist examples for which
\begin{align*}
\hat{\rho}_{\mathcal{V}},_{G_{2}}\left(  \mathcal{A}\right)    &
\neq\hat
{\rho}_{\mathcal{V}},_{G_{2}^{\prime}}\left(  \mathcal{A}\right)  ,\text{ }\\
\hat{\rho}_{\mathcal{V}},_{\overline{G}_2}\left(
\mathcal{A}\right)    & \neq\hat
{\rho}_{\mathcal{V}},_{\overline{G}_2^{\prime}}\left(
\mathcal{A}\right)  .
\end{align*}
\end{remark}

The diagram in Figure~\ref{fig:Hasse.diag} summarizes the results
of this section. We remark that no relations other than the ones
given in Figure~\ref{fig:Hasse.diag} can be established among
these path-complete graphs. Indeed, whenever there are no
relations between two graphs in Figure~\ref{fig:Hasse.diag}, we
have examples of matrices $A_1,A_2$ for which
one graph can outperform the other. \amal{These examples are not presented here but are available online and can be retrieved from \cite{uclurl}.}

Based on our numerical experiments, the graphs
$G_1$ and $G^\prime_1$ seem to statistically perform better than
all other graphs in Figure~\ref{fig:Hasse.diag}. For example, we
ran experiments on a set of $100$ random $5\times5$ matrices
$\{A_1,A_2\}$ with elements uniformly distributed in $\left[
-1,1\right]$ to compare the performance of graphs $G_1, G_2$ and
$\overline{G}_2$. If in each case we also consider the relabeled
matrices (i.e., $\{A_2,A_1\}$) as our input, then, out of the
total $200$ instances, graph $G_1$ produced strictly better bounds
on the JSR $58$ times, whereas graphs $G_2$ and $\overline{G}_2$
each produced the best bound of the three graphs only $23$ times.
(The numbers do not add up to $200$ due to ties.) In addition to
this superior performance, the bound
$\hat{\rho}_{\mathcal{V}},_{G_{1}}\left( \left\{
A_{1},A_{2}\right\}  \right)  $ obtained by analysis via the graph
$G_1$ is invariant under (i) permutation of the labels $A_{1}$ and
$A_{2}$ (obvious), and (ii) transposing of $A_{1}$ and $A_{2}$
(Corollary~\ref{cor:HSCC.invariance.under.transpose}). These are
desirable properties which fail to hold for $G_{2}$ and
$\overline{G}_2$ or their duals. Motivated by these observations,
we generalize $G_{1}$ and its dual $G_{1}^{\prime}$ in the next
section to the case of $m$ matrices and $m$ Lyapunov functions and
establish that they have certain appealing properties. We will
prove (cf. Theorem \ref{thm:HSCC.bound}) that these graphs always
perform better than a common Lyapunov function in 2 steps (i.e.,
the graph $H_2$ in Figure~\ref{fig:jsr.graphs}), whereas, this is
not the case for $G_{2}$ and $\overline{G}_2$ or their duals.


\section{Further analysis of a particular family of path-complete graphs}\label{sec:hscc}
The framework of path-complete graphs provides a multitude of
semidefinite programming based techniques for the approximation of
the JSR whose performance vary with computational cost. For
instance, as we increase the number of nodes of the graph, or the
degree of the polynomial Lyapunov functions assigned to the nodes,
or the number of edges of the graph that instead of labels of
length one have labels of higher length, we obtain better results
but at a higher computational cost. Many of these approximation
techniques are asymptotically tight, so in theory they can be used
to achieve any desired accuracy of approximation. For example,
$$\hat{\rho}_{\mathcal{V}^{SOS,2d}}(\mathcal{A})\rightarrow\rho(\mathcal{A})
\ \mbox{as} \ 2d\rightarrow\infty,$$ where $\mathcal{V}^{SOS,2d}$
denotes the class of sum of squares homogeneous polynomial
Lyapunov functions of degree $2d$. (Recall our notation for bounds
from Section~\ref{subsec:comparison.of.who.beats.who}.) It is also
true that a common quadratic Lyapunov function for products of
higher length achieves the true JSR
asymptotically~\amal{\cite{BlimanFerrari}},~\cite{Raphael_Book}; i.e.\footnote{By
$\mathcal{V}^2$ we denote the class of quadratic homogeneous
polynomials. We drop the superscript ``SOS'' because nonnegative
quadratic polynomials are always sums of squares.},
$$\sqrt[t]{\hat{\rho}_{\mathcal{V}^{2}}(\mathcal{A}^t)}\rightarrow\rho(\mathcal{A})\ \mbox{as} \ t\rightarrow\infty.$$

Nevertheless, it is desirable for practical purposes to identify a
class of path-complete graphs that provide a good tradeoff between
quality of approximation and computational cost. Towards this
objective, we propose the use of $m$ quadratic Lyapunov functions
assigned to the nodes of the De Bruijn graph\footnote{The De
Bruijn graph of dimension $k$ on $m$ symbols is a labeled directed
graph with $m^k$ nodes and $m^{k+1}$ edges whose nodes are indexed
by all possible words of length $k$ from the alphabet $\{1,\ldots,
m\}$, and whose edges have labels of length one and are obtained
by the following simple rule: There is an edge labeled with the
letter $j$ (or for our purposes the matrix $A_j$) going from node
$i_1i_2\ldots i_{k-1}i_k$ to node $i_2i_3\ldots i_{k}j$, $\forall
i_1\ldots i_k\in\{1,\ldots,m\}^k$ and $\forall
j\in\{1,\ldots,m\}$.} of dimension $1$ on $m$ symbols for the
approximation of the JSR of a set of $m$ matrices. \amal{This is precisely the graph of path-dependent quadratic Lyapunov functions of path length $1$~\cite{LeeD06}.} This graph and
its dual are particular path-complete graphs with $m$ nodes and
$m^2$ edges and will be the subject of study in this section. If
we denote the quadratic Lyapunov functions by $x^TP_ix$, then we
are proposing the use of linear matrix inequalities
\begin{equation}\label{eq:HSCC.lmis}
\begin{array}{rll}
P_i&\succ&0 \quad  \forall i=1,\ldots, m, \\
\gamma^2A_i^TP_jA_i&\preceq&P_i \quad \forall i,j=\{1,\ldots,
m\}^2
\end{array}
\end{equation}
or the set of LMIs
\begin{equation}\label{eq:dual.HSCC.lmis}
\begin{array}{rll}
P_i&\succ&0 \quad  \forall i=1,\ldots, m, \\
\gamma^2A_i^TP_iA_i&\preceq&P_j \quad \forall i,j=\{1,\ldots,
m\}^2
\end{array}
\end{equation}
for the approximation of the JSR of $m$ matrices. \amal{We note that the LMIs in (\ref{eq:HSCC.lmis}) have appeared in~\cite{daafouzbernussou,LeeK09} and those in (\ref{eq:dual.HSCC.lmis}) have appeared in~\cite{LeeD06}}. Throughout this section, we denote the path-complete graphs associated with
(\ref{eq:HSCC.lmis}) and (\ref{eq:dual.HSCC.lmis}) with $G_1$ and
$G_1^\prime$ respectively. (The De Bruijn graph of dimension $1$, by
standard convention, is actually the graph $G_1^\prime$.) Observe
that $G_1$ and $G_1^\prime$ are indeed dual graphs as they can be
obtained from each other by reversing the direction of the edges.
For the case $m=2$, our notation is consistent with the previous
section and these graphs are illustrated in
Figure~\ref{fig:dual.graphs}. Also observe from
Corollary~\ref{cor:min.of.quadratics} and
Corollary~\ref{cor:max.of.quadratics} that the LMIs in
(\ref{eq:HSCC.lmis}) give rise to max-of-quadratics Lyapunov
functions, whereas the LMIs in (\ref{eq:dual.HSCC.lmis}) lead to
min-of-quadratics Lyapunov functions. We will prove in this
section that the approximation bound obtained by these LMIs (i.e.,
the reciprocal of the largest $\gamma$ for which the LMIs
(\ref{eq:HSCC.lmis}) or (\ref{eq:dual.HSCC.lmis}) hold) is always
the same and lies within a multiplicative factor of
$\frac{1}{\sqrt[4]{n}}$ of the true JSR, where $n$ is the
dimension of the matrices. The relation between the bound obtained
by a pair of dual path-complete graphs has a connection to
transposition of the matrices in the set $\mathcal{A}$. We explain
this next.

\subsection{Duality and invariance under transposition}\label{subsec:duality.and.transposition}
In~\cite{dual_LMI_diff_inclusions},~\cite{convex_conjugate_Lyap},
it is shown that absolute asymptotic stability of the linear
difference inclusion in (\ref{eq:linear.difference.inclusion})
defined by the matrices $\mathcal{A}=\{A_1,\ldots,A_m\}$ is
equivalent to absolute asymptotic stability of
(\ref{eq:linear.difference.inclusion}) for the transposed matrices
$\mathcal{A}^T\mathrel{\mathop:}=\{A_1^T,\ldots,A_m^T\}$. Note
that this fact is immediately seen from the definition of the JSR
in (\ref{eq:def.jsr}), since
$\rho(\mathcal{A})=\rho(\mathcal{A}^T)$. It is also well-known
that
$$\hat{\rho}_{\mathcal{V}^2}(\mathcal{A})=\hat{\rho}_{\mathcal{V}^2}(\mathcal{A}^T).$$
Indeed, if $x^TPx$ is a common quadratic Lyapunov function for the
set $\mathcal{A}$, then it is easy to show that $x^TP^{-1}x$ is a
common quadratic Lyapunov function for the set $\mathcal{A}^T$.
However, this nice property is not true for the bound obtained
from some other techniques. For instance, the next example shows
that
\begin{equation}\label{eq:bound.sos.not.transp.invar}
\hat{\rho}_{\mathcal{V}^{SOS,4}}(\mathcal{A})\neq\hat{\rho}_{\mathcal{V}^{SOS,4}}(\mathcal{A}^T),
\end{equation}
i.e., the upper bound obtained by searching for a common quartic
SOS polynomial is not invariant under transposition.

\begin{example}\label{ex:quartic.sos.not.transpose.invariant}
Consider the set of matrices $\mathcal{A}=\{A_1,A_2,A_3,A_4\},$
with \scalefont{.5}
\begin{align*}
A_{1} =\left[
\begin{array}
[c]{rrr}%
10 & -6 & -1 \\
8 & 1 & -16 \\
-8 & 0 & 17
\end{array}
\right] , A_{2}=\left[
\begin{array}
[c]{rrr}%
-5 & 9 & -14\\
1 & 5 & 10 \\
3 & 2 & 16
\end{array}
\right],  A_{3} =\left[
\begin{array}
[c]{rrr}%
-14 & 1 & 0\\
-15 & -8 & -12 \\
-1 & -6 & 7
\end{array}
\right],  A_{4} =\left[
\begin{array}
[c]{rrr}%
1 & -8 & -2\\
1 & 16 & 3 \\
16 & 11 & 14
\end{array}
\right].
\end{align*}
\normalsize We have
$\hat{\rho}_{\mathcal{V}^{SOS,4}}(\mathcal{A})=21.411,$ but
$\hat{\rho}_{\mathcal{V}^{SOS,4}}(\mathcal{A}^T)=21.214$ (up to
three significant digits). \amal{This phenomenon is \emph{not} due to the SOS relaxation and should be attributed to common quartic polynomial Lyapunov functions more generally. We know this because all five polynomial nonnegativity conditions in this problem (on the Lyapunov function and its decrements w.r.t. the four matrices) are imposed on ternary quartic forms. It is known from an old result of Hilbert~\cite{Hilbert_1888} that all nonnegative ternary quartic forms are SOS.}
\end{example}

Similarly, the bound obtained by non-convex inequalities proposed
in~\cite{dual_LMI_diff_inclusions} is not invariant under
transposing the matrices. For such methods, one would have to run
the numerical optimization twice---once for the set $\mathcal{A}$
and once for the set $\mathcal{A}^T$---and then pick the better
bound of the two. We will show that by contrast, the bound
obtained from the LMIs in (\ref{eq:HSCC.lmis}) and
(\ref{eq:dual.HSCC.lmis}) are invariant under transposing the
matrices. Before we do that, let us prove a general result which
states that for path-complete graphs with quadratic Lyapunov
functions as nodes, transposing the matrices has the same effect
as dualizing the graph. \amal{We are grateful to a reviewer who kindly made us aware that an independent and earlier proof of this fact for certain families of path-complete graphs appears in~\cite{LeeK09}.}

\begin{theorem}\label{thm:transpose.bound.dual.bound}
Let $G(N,E)$ be a path-complete graph, and let
$G^\prime(N,E^\prime)$ be its dual graph. Then,
\begin{equation}\label{eq:rho.hat.A^T=rho.hat.G'}
\hat{\rho}_{\mathcal{V}^2,
G}(\mathcal{A}^T)=\hat{\rho}_{\mathcal{V}^2,
G^\prime}(\mathcal{A}).
\end{equation}
\end{theorem}

\begin{proof}
For ease of notation, we prove the claim for the case where the
edge labels of $G(N,E)$ have length one. The proof of the general
case is identical. Pick an arbitrary edge $(i,j)\in E$ going from
node $i$ to node $j$ and labeled with some matrix
$A_l\in\mathcal{A}$. By the application of the Schur complement we
have
\begin{equation}\nonumber
A_lP_jA_l^T\preceq P_i \ \Leftrightarrow\
\begin{bmatrix}P_i & A_l\\ A_l^T & P_j^{-1} \end{bmatrix}\succeq0 \ \Leftrightarrow\
A_l^TP_i^{-1}A_l\preceq P_j^{-1}.
\end{equation}
But this already establishes the claim since we see that $P_i$ and
$P_j$ satisfy the LMI associated with edge $(i,j)\in E$ when the
matrix $A_l$ is transposed if and only if $P_j^{-1}$ and
$P_i^{-1}$ satisfy the LMI associated with edge $(j,i)\in
E^\prime$.
%
\end{proof}

\begin{corollary}\label{cor:transpose.eq.iff.dual.eq}
$\hat{\rho}_{\mathcal{V}^2,
G}(\mathcal{A})=\hat{\rho}_{\mathcal{V}^2, G}(\mathcal{A}^T)$ if
and only if $\hat{\rho}_{\mathcal{V}^2,
G}(\mathcal{A})=\hat{\rho}_{\mathcal{V}^2,
G^\prime}(\mathcal{A})$.
\end{corollary}

\begin{proof}
This is an immediate consequence of the equality in
(\ref{eq:rho.hat.A^T=rho.hat.G'}).
\end{proof}

It is an interesting question for future research to characterize
the path-complete graphs for which one has
$\hat{\rho}_{\mathcal{V}^2,
G}(\mathcal{A})=\hat{\rho}_{\mathcal{V}^2, G}(\mathcal{A}^T).$ For
example, the above corollary shows that this is obviously the case
for any path-complete graph that is self-dual. Let us show next
that this is also the case for graphs $G_1$ and $G_1^\prime$
despite the fact that they are not self-dual.

\begin{corollary}\label{cor:HSCC.invariance.under.transpose}
For the path-complete graphs $G_1$ and $G_1^\prime$ associated
with the inequalities in (\ref{eq:HSCC.lmis}) and
(\ref{eq:dual.HSCC.lmis}), and for any class of continuous,
homogeneous, and positive definite functions $\mathcal{V}$, we
have
\begin{equation}\label{eq:bound.G=bound.G'.general.Lyap}
\hat{\rho}_{\mathcal{V},G_1}(\mathcal{A})=\hat{\rho}_{\mathcal{V},G_1^\prime}(\mathcal{A}).
\end{equation}
Moreover, if quadratic Lyapunov functions are assigned to the
nodes of $G_1$ and $G_1^\prime$, then we have
\begin{equation}\label{eq:bound.G1=G1'=G1A^T=G1'A^T}
\hat{\rho}_{\mathcal{V}^2,
G_1}(\mathcal{A})=\hat{\rho}_{\mathcal{V}^2,
G_1}(\mathcal{A}^T)=\hat{\rho}_{\mathcal{V}^2,
G_1^\prime}(\mathcal{A})=\hat{\rho}_{\mathcal{V}^2,
G_1^\prime}(\mathcal{A}^T).
\end{equation}
\end{corollary}

\begin{proof}
The proof of (\ref{eq:bound.G=bound.G'.general.Lyap}) is
established by observing that the GLFs associated with $G_1$ and
$G_1^\prime$ can be derived from one another via
$V_i^\prime(A_ix)=V_i(x).$ (Note that we are relying here on the
assumption that the matrices $A_i$ are invertible, which as we
noted in Remark~\ref{rmk:invertibility}, is not a limiting
assumption.) Since (\ref{eq:bound.G=bound.G'.general.Lyap}) in
particular implies that $\hat{\rho}_{\mathcal{V}^2,
G_1}(\mathcal{A})=\hat{\rho}_{\mathcal{V}^2,
G_1^\prime}(\mathcal{A})$, we get the rest of the equalities in
(\ref{eq:bound.G1=G1'=G1A^T=G1'A^T}) immediately from
Corollary~\ref{cor:transpose.eq.iff.dual.eq} and this finishes the
proof. For concreteness, let us also prove the leftmost equality
in (\ref{eq:bound.G1=G1'=G1A^T=G1'A^T}) directly. Let $P_i$,
$i=1,\ldots,m,$ satisfy the LMIs in (\ref{eq:HSCC.lmis}) for the
set of matrices $\mathcal{A}$. Then, the reader can check that
$$\tilde{P}_i=A_i P_i^{-1}A_i^T,\quad i=1,\ldots,m,$$
satisfy the LMIs in (\ref{eq:HSCC.lmis}) for the set of matrices
$\mathcal{A}^T$.
\end{proof}


\subsection{An approximation guarantee}\label{subsec:HSCC.bound}
The next theorem gives a bound on the quality of approximation of
the estimate resulting from the LMIs in (\ref{eq:HSCC.lmis}) and
(\ref{eq:dual.HSCC.lmis}). Since we have already shown that
$\hat{\rho}_{\mathcal{V}^2,
G_1}(\mathcal{A})=\hat{\rho}_{\mathcal{V}^2,
G_1^\prime}(\mathcal{A}),$ it is enough to prove this bound for
the LMIs in (\ref{eq:HSCC.lmis}).

\begin{theorem}\label{thm:HSCC.bound}
Let $\mathcal{A}$ be a set of $m$ matrices in $\mathbb{R}^{n\times
n}$ with JSR $\rho(\mathcal{A})$. Let $\hat{\rho}_{\mathcal{V}^2,
G_1}(\mathcal{A})$ be the bound on the JSR obtained from the LMIs
in (\ref{eq:HSCC.lmis}). Then,
\begin{equation}\label{eq:HSCC.bound.4throot.of.n}
\frac{1}{\sqrt[4]{n}}\hat{\rho}_{\mathcal{V}^2,
G_1}(\mathcal{A})\leq\rho(\mathcal{A})\leq\hat{\rho}_{\mathcal{V}^2,
G_1}(\mathcal{A}).
\end{equation}
\end{theorem}

\begin{proof}
The right inequality is just a consequence of $G_1$ being a
path-complete graph
(Theorem~\ref{thm:path.complete.implies.stability}). To prove the
left inequality, consider the set $\mathcal{A}^2$ consisting of
all $m^2$ products of length two. In view of (\ref{eq:CQ.bound}),
a common quadratic Lyapunov function for this set satisfies the
bound
\begin{equation}\nonumber
\frac{1}{\sqrt{n}}\hat{\rho}_{\mathcal{V}^2}(\mathcal{A}^2)\leq\rho(\mathcal{A}^2).
\end{equation}
It is easy to show that
$$\rho(\mathcal{A}^2)=\rho^2(\mathcal{A}).$$
See e.g.~\cite{Raphael_Book}. Therefore,
\begin{equation}\label{eq:bound.of.CQ.2.steps}
\frac{1}{\sqrt[4]{n}}\hat{\rho}_{\mathcal{V}^2}^{\frac{1}{2}}(\mathcal{A}^2)\leq\rho(\mathcal{A}).
\end{equation}
Now suppose for some $\gamma>0$, $x^TQx$ is a common quadratic
Lyapunov function for the matrices in $\mathcal{A}_\gamma^2$;
i.e., it satisfies
\begin{equation}\nonumber
\begin{array}{rll}
Q&\succ&0 \\
\gamma^4(A_iA_j)^TQA_iA_j&\preceq&Q \quad \forall i,j=\{1,\ldots,
m\}^2.
\end{array}
\end{equation}
Then, we leave it to the reader to check that
\begin{equation}\nonumber
P_i=Q+A_i^TQA_i,\quad i=1,\ldots,m
\end{equation}
satisfy (\ref{eq:HSCC.lmis}). Hence,
\begin{equation}\nonumber
\hat{\rho}_{\mathcal{V}^2,
G_1}(\mathcal{A})\leq\hat{\rho}_{\mathcal{V}^2}^{\frac{1}{2}}(\mathcal{A}^2),
\end{equation}
and in view of (\ref{eq:bound.of.CQ.2.steps}) the claim is
established.
\end{proof}

Note that the bound in (\ref{eq:HSCC.bound.4throot.of.n}) is
independent of the number of matrices. Moreover, we remark that
this bound is tighter, in terms of its dependence on $n$, than the
known bounds for $\hat{\rho}_{\mathcal{V}^{SOS,2d}}$ for any
finite degree $2d$ of the sum of squares polynomials. The reader
can check that the bound in (\ref{eq:SOS.bound}) goes
asymptotically as $\frac{1}{\sqrt{n}}$. Numerical evidence
suggests that the performance of both the bound obtained by sum of
squares polynomials and the bound obtained by the LMIs in
(\ref{eq:HSCC.lmis}) and (\ref{eq:dual.HSCC.lmis}) is much better
than the provable bounds in (\ref{eq:SOS.bound}) and in
Theorem~\ref{thm:HSCC.bound}. The problem of improving these
bounds or establishing their tightness is open. It goes without
saying that instead of quadratic functions, we can associate sum
of squares polynomials to the nodes of $G_1$ and obtain a more
powerful technique for which we can also prove better bounds with
the exact same arguments.

\subsection{Numerical examples and applications}\label{subsec:numerical.examples}
In the proof of Theorem~\ref{thm:HSCC.bound}, we essentially
showed that the bound obtained from LMIs in (\ref{eq:HSCC.lmis})
is tighter than the bound obtained from a common quadratic applied
to products of length two. Our first example shows that the LMIs
in (\ref{eq:HSCC.lmis}) can in fact do better than a common
quadratic applied to products of \emph{any} finite length. \amal{We remind the reader that these LMIs correspond to the dual of the De Bruijn graph of dimension one and appear in \cite{daafouzbernussou},~\cite{LeeK09}.}


\begin{example}
Consider the set of matrices $\mathcal{A}=\{A_1,A_2\},$ with
\[
A_{1}=\left[
\begin{array}
[c]{cc}%
1 & 0\\
1 & 0
\end{array}
\right]  ,\text{ }A_{2}=\left[
\begin{array}
[c]{cr}%
0 & 1\\
0 & -1
\end{array}
\right].
\]
This is a benchmark set of matrices that has been studied
in~\cite{Ando98},~\cite{Pablo_Jadbabaie_JSR_journal},~\cite{AAA_PP_CDC08_non_monotonic}
because it gives the worst case approximation ratio of a common
quadratic Lyapunov function. Indeed, it is easy to show that
$\rho(\mathcal{A})=1$, but
$\hat{\rho}_{\mathcal{V}^2}(\mathcal{A})=\sqrt{2}$. Moreover, the
bound obtained by a common quadratic function applied to the set
$\mathcal{A}^t$ is
$$\hat{\rho}_{\mathcal{V}^2}^{\frac{1}{t}}(\mathcal{A}^t)=2^{\frac{1}{2t}},$$
which for no finite value of $t$ is exact. On the other hand, we
show that the LMIs in (\ref{eq:HSCC.lmis}) give the exact bound;
i.e., $\hat{\rho}_{\mathcal{V}^2, G_1}(\mathcal{A})=1$. Due to the
simple structure of $A_{1}$ and $A_{2}$, we can even give an
analytical expression for our Lyapunov functions. Given any
$\varepsilon>0$, the LMIs in (\ref{eq:HSCC.lmis}) with
$\gamma=1/\left( 1+\varepsilon \right)  $ are feasible with
\[
P_{1}=\left[
\begin{array}
[c]{cc}%
a & 0\\
0 & b
\end{array}
\right]  ,\text{\qquad}P_{2}=\left[
\begin{array}
[c]{cc}%
b & 0\\
0 & a
\end{array}
\right],
\]
for any $b>0$ and $a>b/2\varepsilon.$
\end{example}

\begin{example}
\aaa{ Consider the set of randomly generated matrices
$\mathcal{A}=\{A_1,A_2,A_3\},$ with \scalefont{.5}
\begin{align*}
A_{1} =\left[
\begin{array}
[c]{rrrrr}%
0 & -2 & 2 & 2 & 4\\
0 & 0 & -4 & -1 & -6\\
2 & 6 & 0 & -8 & 0\\
-2 & -2 & -3 & 1 & -3\\
-1 & -5 & 2 & 6 & -4
\end{array}
\right] , A_{2}=\left[
\begin{array}
[c]{rrrrr}%
-5 & -2 & -4 & \text{ \ \hspace*{0.01in}}6 & -1\\
1 & 1 & 4 & 3 & -5\\
-2 & 3 & -2 & 8 & -1\\
0 & 8 & -6 & 2 & 5\\
-1 & -5 & 1 & 7 & -4
\end{array}
\right],  A_{3} =\left[
\begin{array}
[c]{rrrrr}%
3 & -8 & -3 & 2 & -4\\
-2 & -2 & -9 & 4 & -1\\
2 & 2 & -5 & -8 & 6\\
-4 & -1 & 4 & -3 & 0\\
0 & 5 & 0 & -3 & 5
\end{array}
\right].
\end{align*}
\normalsize

}

A lower bound on $\rho(\mathcal{A})$ is
$\rho(A_{1}A_{2}A_{2})^{1/3}=11.8015$. The upper approximations
for $\rho(\mathcal{A})$ that we computed for this example are
\aaa{as follows}:
\begin{equation}
\begin{array}{rll}
\hat{\rho}_{\mathcal{V}^2}(\mathcal{A})&=& 12.5683   \\
\hat{\rho}_{\mathcal{V}^2}^{\frac{1}{2}}(\mathcal{A}^2)&=&  11.9575    \\
\hat{\rho}_{\mathcal{V}^2,G_1}(\mathcal{A})&=&   11.8097 \\
\hat{\rho}_{\mathcal{V}^{SOS,4}}(\mathcal{A})&=&   11.8015.
\end{array}
\end{equation}
The bound $\hat{\rho}_{\mathcal{V}^{SOS,4}}$ matches the lower
bound numerically and is most likely exact for this example. This
bound is slightly better than $\hat{\rho}_{\mathcal{V}^2,G_1}$.
However, a simple calculation shows that the semidefinite program
resulting in $\hat{\rho}_{\mathcal{V}^{SOS,4}}$ has 25 more
decision variables than the one for
$\hat{\rho}_{\mathcal{V}^2,G_1}$. Also, the running time of the
algorithm leading to $\hat{\rho}_{\mathcal{V}^{SOS,4}}$ is
noticeably larger than the one leading to
$\hat{\rho}_{\mathcal{V}^2,G_1}$. In general, when the dimension
of the matrices is large, it can often be cost-effective to
increase the number of the nodes of our path-complete graphs but
keep the degree of the polynomial Lyapunov functions assigned to
its nodes relatively low. \amal{For example, a path-dependent quadratic Lyapunov function with path length $2$ (i.e. the De Bruijn of dimension $2$) also achieves the exact JSR by solving a system of LMIs with $9$ quadratic functions and $27$ constraints.}
%
\end{example}
%

\amal{
\begin{example}
\aaa{Consider the set of matrices
$\mathcal{A}=\{A_1,A_2\},$ with \scalefont{1}
\begin{align*}
 A_{1}=\left[
\begin{array}
[c]{rrrrr}%
-1 & -1 \\
-4  & 0
\end{array}
\right],  A_{2} =\left[
\begin{array}
[c]{rrrrr}%
3 & 3 \\
-2 & 1 
\end{array}
\right].
\end{align*}
\normalsize
}
A lower bound for $\rho(\mathcal{A})$ is
$\rho(A_{2}A_{1})^{1/2}=3.917384715148$. Here are some upper approximations for this example computed via four methods:
\begin{equation}
\begin{array}{rll}
\hat{\rho}_{\mathcal{V}^2}^{\frac{1}{2}}(\mathcal{A}^2)&=&  3.9264    \\
\hat{\rho}_{\mathcal{V}^{SOS,4}}(\mathcal{A})&=&   3.9241 \\
\hat{\rho}_{\mathcal{V}^2,G_1}(\mathcal{A})&=&   3.9224 \\
\hat{\rho}_{\mathcal{V}^2,H_3}(\mathcal{A})&=&   3.917384715148.
\end{array}
\end{equation}
This example is interesting because the graph $H_3$ (see Fig. \ref{fig:jsr.graphs}) is the cheapest computational method among the four (e.g. it has only one unknown matrix variable and three constraints, versus one unknown and four constraints for $H_2,$ two unknowns and four constraints for $G_1$), but yet it is the only method that gets the JSR exactly. This shows that the quality of the different methods depends on the particular set of matrices.  In particular, the method corresponding to the graph $H_3,$ which has not appeared in the literature to the best of our knowledge, can outperform other choices in many randomly generated examples.
For this example, if we increase the degree of the common SOS Lyapunov function from $4$ to $6$, or the path length of the path-dependent quadratic Lyapunov function from $1$ to $2$, then these methods also get the JSR exactly though at a higher computational cost.
\end{example}}

\mr{
\begin{example}
\aaa{ Consider the set of matrices
$\mathcal{A}=\{A_1,A_2\},$ with \scalefont{1}
\begin{align*}
 A_{1}=\left[
\begin{array}
[c]{rrrrr}%
0.8 & 0.65 \\
-0.34  & 0.9
\end{array}
\right],  A_{2} =\left[
\begin{array}
[c]{rrrrr}%
0.43 & 0.62 \\
-1.48 & 0.14 
\end{array}
\right].
\end{align*}
\normalsize
}

A lower bound for $\rho(\mathcal{A})$ is
$\rho(A_{1}A_{1}A_{1}A_{2})^{1/4}=1.1644.$ \amal{Here are three upper bounds computed for this example:
\begin{equation}
\begin{array}{rll}
\hat{\rho}_{\mathcal{V}^2}^{\frac{1}{2}}(\mathcal{A}^2)&=&  1.2140    \\
\hat{\rho}_{\mathcal{V}^2,G_1}(\mathcal{A})&=&   1.1927\\
\hat{\rho}_{\mathcal{V}^2,H_3}(\mathcal{A})&=&   1.1875. \\
\end{array}
\end{equation}
Once again, graph $H_3$, which is an example of a new method, outperforms the other two methods even though it solves a smaller semidefinite program.

What is also interesting in the above example is that it is quite challenging to prove that $\rho(A_{1}A_{1}A_{1}A_{2})^{1/4}$ in fact gives the exact JSR; i.e. it is hard to find a matching upper bound. This goal can be achieved for example by a common SOS Lyapunov function of degree $14$, but not by one of degree $12$ or lower. Similarly, path dependent quadratic Lyapunov functions of path lengths $1$, $2$, $3$, \amral{or $4$} fail to find the exact JSR. However, if we combine the SOS method with path-dependent Lyapunov functions (i.e. assign SOS Lyapunov functions to nodes of the De Bruijn graph), then the exact JSR can be achieved by ``\{path length, SOS degree\} pairs'' equal to $\{1,10\}$ or $\{2,8\}$ or $\{3,6\}$.

\amral{If one works with quadratic Lyapunov functions only, then path dependent quadratic Lyapunov functions of path length $5$ succeed in getting the JSR exactly. The resulting semidefinite program has $32$ unknown Lyapunov functions (matrix variables) and $96$ LMIs. By using new path-complete graphs, we were able to get the JSR exactly with only $6$ unknown quadratic Lyapunov functions and $42$ LMIs. The graph that achieved this (not shown) consists of $6$ nodes and $36$ edges and is closely related to Remark~\ref{rmk:nontrivial.graph.generalized}. Each node of this graph has $6$ outgoing edges with exactly the same label going to the $6$ nodes of the graph. The labels on the outgoing edges of the different nodes are respectively $\{A_2,A_1A_2,A_1^2A_2,A_1^3A_2,A_1^4A_2,A_1^5\}$. We leave it to the reader to check that this graph is path-complete.} }
\end{example}}

\amral{

{\bf  Performance on application-motivated problems.} In the remainder of this section, we consider computational problems that arise from three different application scenarios. In all of these applications, the underlying problems have been already shown by the existing literature to be related to the computation of the JSR of certain matrices. We thus focus on the computational aspects, and demonstrate the usefulness of the path-complete graph Lyapunov function framework in situations that arise from practical scenarios. 

\begin{figure}
\centering \includegraphics[scale=0.55]{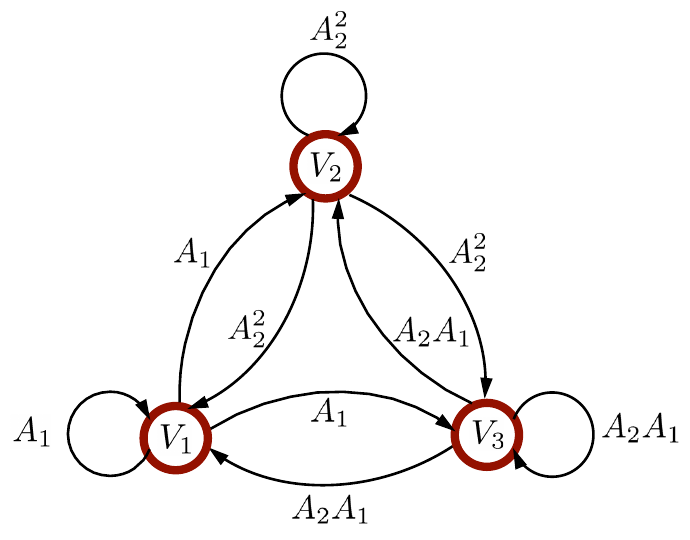}
\caption{\amral{The path-complete graph ${LH}_3$ is obtained from $H_3$ by associating each word in the set $\{A_1,A_2A_1,A_2^2\}$ with a different node on a complete directed graph of order 3, in which all the outgoing edges from every node have the same label.}}
\label{fig:FigureXX}
\end{figure}

Before we proceed, we introduce two new graphs $LH_3$ and $LH_3^2$ which can be verified to be path-complete\footnote{For brevity, we do not provide proofs of path-completeness.}. The first graph, $LH_3$, is shown in Figure~\ref{fig:FigureXX} and is obtained by associating each word in the set $\{A_1,A_2A_1,A_2^2\}$ with a different node on a complete directed graph of order 3, in which all the outgoing edges from every node have the same label. The second graph, $LH_3^2$ (not shown), is a complete directed graph of order 9, and is obtained by applying the same principle to the set\footnote{The words in this set correspond to paths of length two on $H_3$.} $\{A_1^2,~A_2^2A_1,~A_2A1^2,~A_2^4,~A_2A_1A_2^2,~A_1A_2^2,~A_2A_1A_2A_1,~A_2^3A_1,~A_1A_2A_1\}.$ Finally, we will use $D_n$ to denote the De Bruijn graph of dimension $n$. (The number of symbols of the De Bruijn graph will be clear from the context and always equal to the number of matrices whose JSR is under consideration.)

\begin{example}{\bf Application to Computation of the Number of Overlap-Free Words.} It was shown in \cite{JPB2009} that the problem of computation of the smallest exponent of growth of the number of overlap-free or repetition-free binary words (see, e.g. \cite{olfwRef}) reduces to the computation of the JSR of two sparse matrices $A_1$ and $A_2$ in $\mathbb{Z}^{20 \times 20}_{+}$. These relatively large-size matrices are explicitly presented in \cite{JPB2009,GP11}, and are not repeated here in the interest of brevity. More precisely, letting $u_n$ denote the number of overlap-free binary words of length $n$, we have: 
\[
\inf\{r~|~\exists C:~u_n \leq Cn^r\}=\log_2\rho(\{A_1,A_2\}).
\]
It was conjectured in \cite{JPB2009} that $\rho(\{A_1,A_2\})=\sqrt{\rho(A_1A_2)}\approx2.51793404.$ This conjecture was recently proven in \cite{GP11} via a variation of the complex polytope algorithm specialized to non-negative matrices. In Table \ref{tab:table0} we report the results of numerical computation of upper bounds on $\rho(\{A_1,A_2\})$ using various path-complete graphs. The approximate solver times are also reported which correspond to the CPU time of a 2.5 Ghz PC running the solver SeDuMi \cite{Strum1999} on MATLAB. 

\begin{table}[h]
\vspace{-0.1in}
  \begin{center}
  \begin{tabular}{ | r || l | c |}
    \hline
         &    & solver time\\ \hline \hline
    $\hat{\rho}_{\mathcal{V}^2,G_1}(\mathcal{A})=$  & 2.5259 & 0.25 sec\\ \hline
    $\hat{\rho}_{\mathcal{V}^2,H_3}(\mathcal{A})=$  & 2.5223 & 0.16 sec\\ \hline
    $\hat{\rho}_{\mathcal{V}^2,D_2}(\mathcal{A})=$  & 2.51793404 & 1.87 sec\\ \hline
    $\hat{\rho}_{\mathcal{V}^2,LH_3}(\mathcal{A})=$  & 2.51793404 & 1.05 sec\\ \hline
    $\rho(\mathcal{A})=$    & 2.51793404 & \\
    \hline
  \end{tabular}
\end{center}
\caption{}
\label{tab:table0}
\end{table}

The graphs $D_2$ and $LH_3$ indeed provide an exact (up to machine precision) numerical value of the JSR and the running time of the SDP associated with $LH_3$ is only $1$ second. These computations show that the path-complete graph Lyapunov function framework can provide very efficient methods for computation of the JSR in situations of practical and theoretical interest. 



\end{example}

\begin{table}[h]
\vspace{-0.1in}
  \begin{center}
  \begin{tabular}{ | r || c | c | c | c | c | c | c | c | c | c | c |}
    \hline
         &  N=17 & N=19 & solver time\\ \hline \hline
    $\hat{\rho}_{\mathcal{V}^2,H_1}(\mathcal{A})=$ & 0.118781760 & 0.097472458 & 0.15 sec\\ \hline
    $\hat{\rho}_{\mathcal{V}^2,H_4}(\mathcal{A})=$ &     & 0.097471788 & 0.36 sec\\ \hline
    $\hat{\rho}_{\mathcal{V}^2,G_1}(\mathcal{A})=$ &     & 0.097463499 & 0.46 sec\\ \hline
    $\hat{\rho}_{\mathcal{V}^2,H_3}(\mathcal{A})=$ &     &  0.097407530& 0.37 sec\\ \hline
    $\hat{\rho}_{\mathcal{V}^2,H_1^3}(\mathcal{A})=$ &  & 0.097403543 & 0.77 sec\\ \hline
    $\hat{\rho}_{\mathcal{V}^2,D_2}(\mathcal{A})=$  &     & 0.097334910 & 0.75 sec\\ \hline
    $\hat{\rho}_{\mathcal{V}^2,D_3}(\mathcal{A})=$  &     & 0.097332287 & 1.12 sec\\ \hline
       $\hat{\rho}_{\mathcal{V}^2,D_6}(\mathcal{A})=$  &   & 0.097306933 & 9.66 sec\\ \hline
    $\hat{\rho}_{\mathcal{V}^2,LH_3}(\mathcal{A})=$ &     & 0.097380084 & 0.60 sec\\ \hline
    $\hat{\rho}_{\mathcal{V}^2,LH^2_3}(\mathcal{A})=$ &  & 0.097306828& 3.70 sec\\ \hline
    $\rho(\mathcal{A})=$   & 0.118781760  & 0.097301716\footnote{This number is only a \emph{lower bound} on the JSR, given by $\rho(A_1^2A_2^2)^{1/4}$. We conjecture that it is equal to the true JSR.} & \\
    \hline
  \end{tabular}
\end{center}
\caption{}
\label{tab:table1}
\end{table}

\begin{example}{\bf Application to Computation of the Euler Ternary Partition Function.}
The problem of computation of the smallest exponent of growth of the Euler ternary partition function \cite{protasov2000} can be reduced to the problem of computation of the JSR of three matrices with binary 0 or 1 entries. Herein, we examine a special case reported in \cite{GP11}, where the complex polytope method is applied to provide the exact value of the JSR of three 7-by-7 matrices with 0 and 1 entries: 
\[
\rho(\{A_1,A_2,A_3\})=\sqrt{\rho(A_2A_3)}\approx4.722045134.
\]
In this case, the path-complete DeBruijn graph of dimension 1 yields an upper bound on the JSR with great accuracy in a fraction of a second; we have $\hat{\rho}_{\mathcal{V}^2,D_1}(\mathcal{A})= 4.722045134$, and the computation time is 0.15sec on a 2.5 Ghz PC.
\end{example}

\begin{example}{\bf Application to Continuity of Wavelet Functions.} Daubechies' wavelet functions are orthonormal functions $\phi_N$ with compact support on $[0,N],$ satisfying 
\[
\phi_N(x)=\sum_{k=0}^Nc_k\phi_N(2x-k)
\]
where, $N$ is a positive integer and the coefficients $c_k,~0\leq k\leq N$, satisfy certain additional constraints \cite{Gripenberg, Raphael_Book}. The problem of computation of the H\"older exponent of continuity of the wavelet functions \cite{DaLa1992two} is closely related to the problem of computation of the JSR of two linear operators, see, e.g., \cite[Chap. 5]{Raphael_Book},~\cite{Wavelet_book_JSR} and \cite{Gripenberg}. Herein, we are interested in computation of the JSR of the associated matrices for values of $N=5,7,\dots,19.$ The matrix pairs $\{A_{1N},A_{2N}\}$ are of dimension $(N-1)/2$ and have been posted online in \cite{uclurl_wavelet} along with annotated MATLAB code for their computation. We remark that the JSR of the associated pairs of matrices for odd values of $N \in [5,15]$, were first reported in \cite{Gripenberg}, where it was shown (numerically) that 
\begin{equation}\label{eqn:r=rmax}
\rho(\{A_{1N},A_{2N}\})=\max(\rho(A_{1N}),\rho(A_{2N})), \quad N=3,\dots,15.
\end{equation} 
Our numerical analysis conforms with the results of \cite{Gripenberg} for $N\leq 15$ and a single common quadratic Lyapunov function ($\hat{\rho}_{\mathcal{V}^2,H_1}(\cdot)$) provides the exact (up to machine precision) numerical value of the JSR. For brevity, we do not repeat here the numerical values of the JSR for $N \leq 15$, and instead present the numerical upper-bound on the JSR for two more values of $N$, i.e., $N=17$ and $N=19$. Table \ref{tab:table1} summarizes\footnote{As before, the approximate solver times correspond to the CPU time of a 2.5 Ghz PC running SeDuMi \cite{Strum1999} on MATLAB. } our numerical analysis for $N=17$ and $N=19$. For $N=17$ the pattern holds and a single common quadratic Lyapunov function provides the exact value of the JSR which also satisfies \eqref{eqn:r=rmax}. Surprisingly, however, for $N=19$ this pattern breaks and not only \eqref{eqn:r=rmax} does not hold, but also a common quadratic Lyapunov function does not give the exact upper bound! The best upper bound we are providing is obtained by graph $LH_3^2$ which has 9 nodes and 90 LMIs. To the best of our knowledge none of the methods in the existing literature provide a better upper bound at a comparable computation cost.

\end{example}

}

\section{Converse Lyapunov theorems and approximation with arbitrary accuracy}\label{sec:converse.thms}
It is well-known that existence of a Lyapunov function which is
the pointwise maximum of quadratics is not only sufficient but
also necessary for absolute asymptotic stability of
(\ref{eq:switched.linear.sys}) or
(\ref{eq:linear.difference.inclusion}); see
e.g.~\cite{pyatnitskiy_max_lyap}.
This is perhaps an intuitive fact if we recall that switched
systems of type (\ref{eq:switched.linear.sys}) and
(\ref{eq:linear.difference.inclusion}) always admit a convex
Lyapunov function. Indeed, if we take ``enough'' quadratics, the
convex and compact unit sublevel set of a convex Lyapunov function
can be approximated arbitrarily well with sublevel sets of
max-of-quadratics Lyapunov functions, which are intersections of
ellipsoids. This of course implies that the bound obtained from
max-of-quadratics Lyapunov functions is asymptotically tight for
the approximation of the JSR. However, this converse Lyapunov
theorem does not answer two natural questions of importance in
practice: (i) How many quadratic functions do we need to achieve a
desired quality of approximation? (ii) Can we search for these
quadratic functions via semidefinite programming or do we need to
resort to non-convex formulations? \amal{The same questions can naturally be asked for min-of-quadratics Lyapunov functions. The theorem and remark that follow provide an answer to these questions by relying on the connections that we have already established between min/max-quadratics Lyapunov functions and path-dependent Lyapunov functions~\cite{LeeD06} and their duals~\cite{LeeK09}. Our results further provides a worst case approximation guarantee for path-dependent quadratic Lyapunov functions of any given path length, and similarly for their duals.}

\begin{theorem}\label{thm:converse.max.of.quadratics}
Let $\mathcal{A}$ be a set of $m$ matrices in $\mathbb{R}^{n\times
n}$. Given any positive integer $l$, there exists an explicit
path-complete graph $G$ consisting of $m^{l-1}$ nodes assigned to
quadratic Lyapunov functions and $m^l$ edges with labels of length
one such that the linear matrix inequalities associated with $G$
imply existence of a max-of-quadratics Lyapunov function and the
resulting bound obtained from the LMIs satisfies
\begin{equation}\label{eq:bound.converse.thm.max}
\frac{1}{\sqrt[2l]{n}}\hat{\rho}_{\mathcal{V}^2,
G}(\mathcal{A})\leq\rho(\mathcal{A})\leq\hat{\rho}_{\mathcal{V}^2,
G}(\mathcal{A}).
\end{equation}
\end{theorem}
\begin{proof}
Let us denote the $m^{l-1}$ quadratic Lyapunov functions by
$x^TP_{i_1\ldots i_{l-1}}x$, where $i_1\ldots
i_{l-1}\in\{1,\ldots,m\}^{l-1}$ is a multi-index used for ease of
reference to our Lyapunov functions. We claim that we can let $G$
be the graph dual to the De Bruijn graph of dimension $l-1$ on $m$
symbols. The LMIs associated to this graph are given by
\begin{equation}\label{eq:converse.thm.LMIs}
\begin{array}{rll}
P_{i_1i_2\ldots i_{l-2}i_{l-1}}&\succ& 0\ \ \ \  \forall i_1\ldots i_{l-1}\in\{1,\ldots,m\}^{l-1}\\
A_j^TP_{i_1i_2\ldots i_{l-2}i_{l-1}}A_j&\preceq&P_{i_2i_3\ldots
i_{l-1}j}\\
\ &\ & \forall i_1\ldots i_{l-1}\in\{1,\ldots,m\}^{l-1},\\
\ &\ & \forall j\in\{1,\ldots,m\}.
\end{array}
\end{equation}
\amal{These LMIs appear in~\cite{LeeK09} and are known to be asymptotically exact.} The fact that $G$ is path-complete and that the LMIs imply
existence of a max-of-quadratics Lyapunov function follows from
Corollary~\ref{cor:max.of.quadratics}. The proof that these LMIs
satisfy the bound in (\ref{eq:bound.converse.thm.max}) is a
straightforward generalization of the proof of
Theorem~\ref{thm:HSCC.bound}. By the same arguments we have
\begin{equation}\label{eq:bound.of.CQ.l.steps}
\frac{1}{\sqrt[2l]{n}}\hat{\rho}_{\mathcal{V}^2}^{\frac{1}{l}}(\mathcal{A}^l)\leq\rho(\mathcal{A}).
\end{equation}
Suppose $x^TQx$ is a common quadratic Lyapunov function for the
matrices in $\mathcal{A}^l$; i.e., it satisfies
\begin{equation}\nonumber
\begin{array}{rll}
Q&\succ&0 \\
(A_{i_1}\ldots A_{i_l})^TQA_{i_1}\ldots A_{i_l}&\preceq&Q \quad
\forall i_1\ldots i_{l}\in\{1,\ldots,m\}^{l}.
\end{array}
\end{equation}
Then, it is easy to check that\footnote{The construction of the
Lyapunov function here is a special case of a general scheme for
constructing Lyapunov functions that are monotonically decreasing
from those that decrease only every few steps; see~\cite[p.
58]{AAA_MS_Thesis}.}
\begin{equation}\nonumber
\begin{array}{lll}
P_{i_1i_2\ldots i_{l-2}i_{l-1}}=Q+A_{i_{l-1}}^TQA_{i_{l-1}}&\ &\ \\
+(A_{i_{l-2}}A_{i_{l-1}})^TQ(A_{i_{l-2}}A_{i_{l-1}})+\cdots&\ &\ \\
+ (A_{i_1}A_{i_2}\ldots
A_{i_{l-2}}A_{i_{l-1}})^TQ(A_{i_1}A_{i_2}\ldots
A_{i_{l-2}}A_{i_{l-1}}), &\ &\ \\
i_1\ldots i_{l-1}\in\{1,\ldots,m\}^{l-1},&\ &\
\end{array}
\end{equation}
satisfy (\ref{eq:converse.thm.LMIs}). Hence,
\begin{equation}\nonumber
\hat{\rho}_{\mathcal{V}^2,
G}(\mathcal{A})\leq\hat{\rho}_{\mathcal{V}^2}^{\frac{1}{l}}(\mathcal{A}^l),
\end{equation}
and in view of (\ref{eq:bound.of.CQ.l.steps}) the claim is
established.
\end{proof}

%

\begin{remark}
\amal{Arbitrarily good approximation bounds identical to those in Theorem~\ref{thm:converse.max.of.quadratics} can be proven for min-of-quadratics Lyapunov functions in a similar fashion}. The only difference is that
the LMIs in (\ref{eq:converse.thm.LMIs}) would get replaced by the
ones corresponding to the dual graph of $G$, \amal{i.e., the De Bruijn graph which is associated with path-dependent Lyapunov functions~\cite{LeeD06}.}
\end{remark}

Our last theorem establishes approximation bounds for a family of
path-complete graphs with one single node but several edges
labeled with words of different lengths. Examples of such
path-complete graphs include graph $H_3$ in
Figure~\ref{fig:jsr.graphs} and graph $H_4$ in
Figure~\ref{fig:non.trivial.path.complete}.

\begin{theorem}
\label{thm-bound-codes} Let $\mathcal{A}$ be a set of matrices in
$\mathbb{R}^{n\times n}.$ Let $\tilde{G}\left(  \left\{ 1\right\}
,E\right)  $ be a path-complete graph, and $l$ be the length of
the shortest word in $\tilde{\mathcal{A}}=\left\{ L\left( e\right)
:e\in E\right\}.$ Then
$\hat{\rho}_{\mathcal{V}^{2}},_{\tilde{G}}(\mathcal{A})$ provides
an estimate of $\rho\left(  \mathcal{A}\right)  $ that satisfies
\[
\frac{1}{\sqrt[2l]{n}}\hat{\rho}_{\mathcal{V}^{2}},_{\tilde{G}}(\mathcal{A}%
)\leq\rho(\mathcal{A})\leq\hat{\rho}_{\mathcal{V}^{2}},_{\tilde{G}%
}(\mathcal{A}).
\]

\end{theorem}

\begin{proof}
The right inequality is obvious, we prove the left one. Since both
$\hat{\rho }_{\mathcal{V}^{2}},_{\tilde{G}}(\mathcal{A})$ and
$\rho$ are homogeneous in $\mathcal{A},$ we may assume, without
loss of generality, that $\hat{\rho
}_{\mathcal{V}^{2}},_{\tilde{G}}(\mathcal{A})=1$. Suppose for the
sake of contradiction that
\begin{equation}\label{eq:supposition-thm-words}\rho(\mathcal{A})<1/\sqrt[2l]{n}.\end{equation}
We will show that this implies that
$\hat{\rho}_{\mathcal{V}^{2}},_{\tilde{G}}(\mathcal{A})<1$.
Towards this goal, let us first prove that
$\rho(\tilde{\mathcal{A}})\leq\rho^l(\mathcal{A}).$
\aaa{ Indeed, if we had
$\rho(\tilde{\mathcal{A}})>\rho^{l}(\mathcal{A})$, then there
would exist\footnote{Here, we are appealing to the well-known fact
about the JSR of a general set of matrices $\mathcal{B}$:
$\rho(\mathcal{B})=\limsup_{k\rightarrow\infty}
\max_{B\in\mathcal{B}^k} \rho^\frac{1}{k}(B).$ See
e.g.~\cite[Chap. 1]{Raphael_Book}.} an integer $i$ and a product
$A_\sigma\in\tilde{\mathcal{A}}^i$ such that
\begin{equation}\label{eq:rho^1/i>rho^l}
\rho^{\frac{1}{i}}(A_\sigma)>\rho^{l}(\mathcal{A}).
\end{equation}
Since we also have $A_\sigma\in\mathcal{A}^j$ (for some $j\geq
il$), it follows that
\begin{equation}\label{eq:rho^1/j<rho}
\rho^{\frac{1}{j}}(A_\sigma)\leq\rho(\mathcal{A}).
\end{equation}
The inequality in (\ref{eq:rho^1/i>rho^l}) together with
$\rho(\mathcal{A})\leq 1$ gives
$$\rho^{\frac{1}{j}}(A_\sigma)>\rho^{\frac{il}{j}}(\mathcal{A})\geq\rho(\mathcal{A}).$$
But this contradicts (\ref{eq:rho^1/j<rho}). Hence we have shown
$$\rho(\tilde{\mathcal{A}})\leq \rho^{l}(\mathcal{A}).$$
}
\rmj{Now, by our hypothesis (\ref{eq:supposition-thm-words})
above, we have that $ \rho(\tilde{\mathcal{A}})<1/\sqrt{n}.$}
Therefore, there
exists $\epsilon>0$ such that $\rho((1+\epsilon)\tilde{\mathcal{A}}%
)<1/\sqrt{n}.$ It then follows from (\ref{eq:CQ.bound}) that there
exists a common quadratic Lyapunov function for
$(1+\epsilon)\tilde{\mathcal{A}}.$ \aaa{Hence}, $\hat
{\rho}_{\mathcal{V}^{2}}((1+\epsilon)\tilde{\mathcal{A}})\leq1,$
which
immediately implies that $\hat{\rho}_{\mathcal{V}^{2}},_{\tilde{G}%
}(\mathcal{A})<1,$ a contradiction.
\end{proof}

\aaa{ A noteworthy immediate corollary of
Theorem~\ref{thm-bound-codes} (obtained by setting
$\tilde{\mathcal{A}}=\bigcup_{t=r}^k \mathcal{A}^t)$ is the
following: If $\rho(\mathcal{A})<\frac{1}{\sqrt[2r]{n}}$, then
there exists a quadratic Lyapunov function that decreases
simultaneously for all products of lengths $r,r+1,\ldots,r+k$, for
any desired value of $k$. Note that this fact is obvious for
$r=1$, but nonobvious for $r\geq 2$.}

\section{Conclusions and future directions}\label{sec:conclusions.future.directions}

We introduced the framework of path-complete graph Lyapunov
functions for the formulation of semidefinite programming based
algorithms for approximating the joint spectral radius (or
equivalently establishing absolute asymptotic stability of an
arbitrarily switched linear system). We defined the notion of a
path-complete graph, which was inspired by concepts in automata
theory. We showed that every path-complete graph gives rise to a
technique for the approximation of the JSR. This provided a
unifying framework that includes many of the previously proposed
techniques and also introduces new ones. \rmj{\aaa{(In fact}, all
families of LMIs that we are aware of \aaa{are} particular cases
of our method.\aaa{)}} \aaa{We shall also emphasize that although
we focused on switched \emph{linear} systems because of our
interest in the JSR, the analysis technique of multiple Lyapunov
functions on path-complete graphs is clearly valid for switched
\emph{nonlinear} systems as well.}

We compared the quality of the bound obtained from certain classes
of path-complete graphs, including all path-complete graphs with
two nodes on an alphabet of two matrices, and also a certain
family of dual path-complete graphs. \amal{Among the different path-complete graphs considered in this paper, we observed that the De Bruijn graph and its dual, whose LMIs appear in the earlier work~\cite{LeeD06},~\cite{LeeK09}, have a superior performance on average (but not always). Motivated by this fact, we studied these graphs in further detail. For example, we showed that stability analysis via these graphs is invariant under transposition of the matrices, results in common min/max-of-quadratics Lyapunov functions, and produces upper bounds on the JSR that are always within a
multiplicative factor of $1/\sqrt[4]{n}$ of the true value, already for the first level of the hierarchy.}
%
%
Finally, we presented two converse Lyapunov theorems, one for the
well-known methods of minimum and maximum-of-quadratics Lyapunov
functions, and the other for a new class of methods that propose
the use of a common quadratic Lyapunov function for a set of words
of possibly different lengths.

We believe the methodology proposed in this \amal{paper} should
straightforwardly extend to the case of \emph{constrained
switching} by requiring the graphs to have a path not for all the
words, but only the words allowed by the constraints on the
switching. A rigorous treatment of this idea is left for future
work.

Another question for future research is to determine the
complexity of checking path-completeness of a given graph
$G(N,E)$. As we explained in Section~\ref{sec:graphs.jsr},
well-known algorithms in automata theory (see e.g.~\cite[Chap.
4]{Hopcroft_Motwani_Ullman_automata_Book}) can check for
path-completeness by testing whether the associated finite
automaton accepts all finite words. When the automata are
deterministic (i.e., when all outgoing edges from every node have
different labels), these algorithms are very efficient and have
running time of only $O(|N|^2)$. However, the problem of deciding
whether a non-deterministic finite automaton accepts all finite
words is known to be PSPACE-complete~\cite[p.
265]{GareyJohnson_Book}. 
%
\amal{O}f course,
the step of checking path-completeness of a graph is done offline
and prior to the run of our algorithms for approximating the JSR.
Therefore, while checking path-completeness is in general
difficult, the approximation algorithms that we presented indeed
run in polynomial time since they work with a fixed (a priori
chosen) path-complete graph. Nevertheless, the question on
complexity of checking path-completeness is interesting in many
other settings, e.g., when deciding whether a given set of
Lyapunov inequalities imply stability of an arbitrarily switched
system.


Some other interesting questions that can be explored in the
future are the following. What are some other classes of
path-complete graphs that lead to new techniques for proving
stability of switched systems? \amral{Can we classify graph operations that preserve path-completeness?} How can we compare the performance
of different path-complete graphs in a systematic way? Given a set
of matrices, a class of Lyapunov functions, and a fixed size for
the graph, can we efficiently come up with the least conservative
topology of a path-complete graph? \rmjj{What properties of a set of matrices make a particular path-complete graph Lyapunov function better than another one?}
\amal{W}hat are the analogues of the
results of this \amal{paper} for continuous time switched systems? To
what extent do the results carry over to the synthesis (controller
design) problem for switched systems? These questions and several
others show potential for much follow-up work on path-complete
graph Lyapunov functions.

\bibliographystyle{abbrv}
\bibliography{pablo_amirali}

\def\cprime{$'$}
\begin{thebibliography}{10}

\bibitem{uclurl}
URL:
  \url{http://perso.uclouvain.be/raphael.jungers/contents/LMI_comparisons.txt}.

\bibitem{uclurl_wavelet}
URL: \url{http://perso.uclouvain.be/raphael.jungers/contents/wavelets.zip}.

\bibitem{AAA_MS_Thesis}
A.~A. Ahmadi.
\newblock Non-monotonic {L}yapunov functions for stability of nonlinear and
  switched systems: theory and computation.
\newblock Master's Thesis, Massachusetts Institute of Technology, June 2008.
  Available from \texttt{http://dspace.mit.edu/handle/1721.1/44206}.

\bibitem{HSCC_JSR_Path_complete}
A.~A. Ahmadi, R.~Jungers, P.~A. Parrilo, and M.~Roozbehani.
\newblock Analysis of the joint spectral radius via {L}yapunov functions on
  path-complete graphs.
\newblock In {\em Hybrid Systems: Computation and Control 2011}, Lecture Notes
  in Computer Science. Springer, 2011.

\bibitem{Path_complete_converse}
A.~A. Ahmadi, R.~M. Jungers, P.~A. Parrilo, and M.~Roozbehani.
\newblock When is a set of {LMI}s a sufficient condition for stability?
\newblock {\em arXiv preprint arXiv:1201.3227}, 2012.

\bibitem{AAA_PP_CDC08_non_monotonic}
A.~A. Ahmadi and P.~A. Parrilo.
\newblock Non-monotonic {L}yapunov functions for stability of discrete time
  nonlinear and switched systems.
\newblock In {\em Proceedings of the 47$^{th}$ IEEE Conference on Decision and
  Control}, 2008.

\bibitem{Ando98}
T.~Ando and M.-H. Shih.
\newblock Simultaneous contractibility.
\newblock {\em SIAM Journal on Matrix Analysis and Applications}, 19:487--498,
  1998.

\bibitem{olfwRef}
J.~Berstel.
\newblock Growth of repetition-free wordsÑa review.
\newblock {\em Theoretical Computer Science}, 340(2):280--290, 2005.

\bibitem{BlimanFerrari}
P.~Bliman and G.~Ferrari-Trecate.
\newblock Stability analysis of discrete-time switched systems through
  {L}yapunov functions with nonminimal state.
\newblock In {\em Proceedings of IFAC Conference on the Analysis and Design of
  Hybrid Systems}, pages 325--330, 2003.

\bibitem{BlNes05}
V.~D. Blondel and Y.~Nesterov.
\newblock Computationally efficient approximations of the joint spectral
  radius.
\newblock {\em SIAM J. Matrix Anal. Appl.}, 27(1):256--272, 2005.

\bibitem{BlNT04}
V.~D. Blondel, Y.~Nesterov, and J.~Theys.
\newblock On the accuracy of the ellipsoidal norm approximation of the joint
  spectral radius.
\newblock {\em Linear Algebra and its Applications}, 394:91--107, 2005.

\bibitem{BlTi2}
V.~D. Blondel and J.~N. Tsitsiklis.
\newblock The boundedness of all products of a pair of matrices is undecidable.
\newblock {\em Systems and Control Letters}, 41:135--140, 2000.

\bibitem{multiple_lyap_Branicky}
M.~S. Branicky.
\newblock Multiple {L}yapunov functions and other analysis tools for switched
  and hybrid systems.
\newblock {\em IEEE Transactions on Automatic Control}, 43(4):475--482, 1998.

\bibitem{Experiment_JSR}
C.~T. Chang and V.~D. Blondel.
\newblock An experimental study of approximation algorithms for the joint
  spectral radius.
\newblock {\em Numerical Algorithms}, pages 1--22, 2012.

\bibitem{daafouzbernussou}
J.~Daafouz and J.~Bernussou.
\newblock Parameter dependent {L}yapunov functions for discrete time systems
  with time varying parametric uncertainties.
\newblock {\em Systems and Control Letters}, 43(5):355--359, 2001.

\bibitem{DaLa1992two}
I.~Daubechies and J.~C. Lagarias.
\newblock Two-scale difference equations ii. local regularity, infinite
  products of matrices and fractals.
\newblock {\em SIAM Journal on Mathematical Analysis}, 23(4):1031--1079, 1992.

\bibitem{jsr_toolbox}
J.~M.~H. G.~Vankeerberghen and R.~M. Jungers.
\newblock The {JSR} toolbox.
\newblock {\em Matlab Central,
  http://www.mathworks.com/matlabcentral/fileexchange/33202-the-jsr-toolbox}.

\bibitem{GareyJohnson_Book}
M.~R. Garey and D.~S. Johnson.
\newblock {\em Computers and Intractability}.
\newblock W. H. Freeman and Co., San Francisco, Calif., 1979.

\bibitem{dual_LMI_diff_inclusions}
R.~Goebel, T.~Hu, and A.~R. Teel.
\newblock Dual matrix inequalities in stability and performance analysis of
  linear differential/difference inclusions.
\newblock In {\em Current Trends in Nonlinear Systems and Control}, pages
  103--122. 2006.

\bibitem{convex_conjugate_Lyap}
R.~Goebel, A.~R. Teel, T.~Hu, and Z.~Lin.
\newblock Conjugate convex {L}yapunov functions for dual linear differential
  inclusions.
\newblock {\em IEEE Transactions on Automatic Control}, 51(4):661--666, 2006.

\bibitem{grip}
G.~Gripenberg.
\newblock Computing the joint spectral radius.
\newblock {\em Linear Algebra and its Applications}, 234:43--60, 1996.

\bibitem{Gripenberg}
G.~Gripenberg.
\newblock Computing the joint spectral radius.
\newblock {\em Linear Algebra and its Applications}, 234:43--60, 1996.

\bibitem{GraphTheory_Handbook}
J.~L. Gross and J.~Yellen.
\newblock {\em Handbook of Graph Theory (Discrete Mathematics and Its
  Applications)}.
\newblock CRC Press, 2003.

\bibitem{GP11}
N.~Guglielmi and V.~Protasov.
\newblock Exact computation of joint spectral characteristics of linear
  operators.
\newblock {\em Foundations of Computational Mathematics}, 13(1):37--97, 2013.

\bibitem{GZalgorithm}
N.~Guglielmi and M.~Zennaro.
\newblock An algorithm for finding extremal polytope norms of matrix families.
\newblock {\em Linear Algebra and its Applications}, 428:2265--2282, 2008.

\bibitem{GuglielmiZennaro2}
N.~Guglielmi and M.~Zennaro.
\newblock Finding extremal complex polytope norms for families of real
  matrices.
\newblock {\em SIAM Journal on Matrix Analysis and Applications},
  31(2):602--620, 2009.

\bibitem{Hilbert_1888}
D.~Hilbert.
\newblock \"{U}ber die {D}arstellung {D}efiniter {F}ormen als {S}umme von
  {F}ormenquadraten.
\newblock {\em Math. Ann.}, 32, 1888.

\bibitem{Hopcroft_Motwani_Ullman_automata_Book}
J.~E. Hopcroft, R.~Motwani, and J.~D. Ullman.
\newblock {\em Introduction to Automata Theory, Languages, and Computation}.
\newblock Addison Wesley, 2001.

\bibitem{composite_Lyap2}
T.~Hu and Z.~Lin.
\newblock Absolute stability analysis of discrete-time systems with composite
  quadratic {L}yapunov functions.
\newblock {\em IEEE Transactions on Automatic Control}, 50(6):781--797, 2005.

\bibitem{composite_Lyap}
T.~Hu, L.~Ma, and Z.~Li.
\newblock On several composite quadratic {L}yapunov functions for switched
  systems.
\newblock In {\em Proceedings of the 45$^{th}$ IEEE Conference on Decision and
  Control}, 2006.

\bibitem{hu-ma-lin}
T.~Hu, L.~Ma, and Z.~Lin.
\newblock Stabilization of switched systems via composite quadratic functions.
\newblock {\em IEEE Transactions on Automatic Control}, 53(11):2571 -- 2585,
  2008.

\bibitem{JohRan_PWQ}
M.~Johansson and A.~Rantzer.
\newblock Computation of piecewise quadratic {L}yapunov functions for hybrid
  systems.
\newblock {\em IEEE Transactions on Automatic Control}, 43(4):555--559, 1998.

\bibitem{Raphael_Book}
R.~Jungers.
\newblock {\em The joint spectral radius: theory and applications}, volume 385
  of {\em Lecture Notes in Control and Information Sciences}.
\newblock Springer, 2009.

\bibitem{JPB2009}
R.~M. Jungers, V.~Y. Protasov, and V.~D. Blondel.
\newblock Overlap-free words and spectra of matrices.
\newblock {\em Theoretical Computer Science}, 410(38):3670--3684, 2009.

\bibitem{LeeD06_disturb}
J.~W. Lee and G.~E. Dullerud.
\newblock Optimal disturbance attenuation for discrete-time switched and
  {M}arkovian jump linear systems.
\newblock {\em SIAM Journal on Control and Optimization}, 45(4):1329--1358,
  2006.

\bibitem{LeeD06}
J.~W. Lee and G.~E. Dullerud.
\newblock Uniform stabilization of discrete-time switched and {M}arkovian jump
  linear systems.
\newblock {\em Automatica}, 42(2):205--218, 2006.

\bibitem{LeeK08_output}
J.~W. Lee and P.~P. Khargonekar.
\newblock Optimal output regulation for discrete-time switched and {M}arkovian
  jump linear systems.
\newblock {\em SIAM Journal on Control and Optimization}, 47(1):40--72, 2008.

\bibitem{LeeK09}
J.~W. Lee and P.~P. Khargonekar.
\newblock Detectability and stabilizability of discrete-time switched linear
  systems.
\newblock {\em IEEE Transactions on Automatic Control}, 54(3):424--437, 2009.

\bibitem{switched_system_survey}
H.~Lin and P.~J. Antsaklis.
\newblock Stability and stabilizability of switched linear systems: a short
  survey of recent results.
\newblock In {\em Proceedings of IEEE International Symposium on Intelligent
  Control}, 2005.

\bibitem{Lind_Marcus_symbolic_Book}
D.~Lind and B.~Marcus.
\newblock {\em An Introduction to Symbolic Dynamics and Coding}.
\newblock Cambridge University Press, 1995.

\bibitem{pyatnitskiy_max_lyap}
A.~Molchanov and Y.~Pyatnitskiy.
\newblock Criteria of asymptotic stability of differential and difference
  inclusions encountered in control theory.
\newblock {\em Systems and Control Letters}, 13:59--64, 1989.

\bibitem{Wavelet_book_JSR}
I.~I. Novikov, V.~I. Protasov, and M.~A. Skopina.
\newblock {\em Wavelet Theory}, volume 239.
\newblock AMS Bookstore, 2011.

\bibitem{Pablo_Jadbabaie_JSR_journal}
P.~A. Parrilo and A.~Jadbabaie.
\newblock Approximation of the joint spectral radius using sum of squares.
\newblock {\em Linear Algebra and its Applications}, 428(10):2385--2402, 2008.

\bibitem{protasov2000}
V.~Y. Protasov.
\newblock Asymptotic behaviour of the partition function.
\newblock {\em Sbornik: Mathematics}, 191(3):381, 2000.

\bibitem{protasov-jungers-blondel09}
V.~Y. Protasov, R.~M. Jungers, and V.~D. Blondel.
\newblock Joint spectral characteristics of matrices: a conic programming
  approach.
\newblock {\em SIAM Journal on Matrix Analysis and Applications},
  31(4):2146--2162, 2010.

\bibitem{MardavijRoozbehani2008}
M.~Roozbehani.
\newblock {\em Optimization of {L}yapunov invariants in analysis and
  implementation of safety-critical software systems}.
\newblock PhD thesis, Massachusetts Institute of Technology, 2008.

\bibitem{Roozbehani2008}
M.~Roozbehani, A.~Megretski, E.~Frazzoli, and E.~Feron.
\newblock Distributed {L}yapunov functions in analysis of graph models of
  software.
\newblock {\em Springer Lecture Notes in Computer Science}, 4981:443--456,
  2008.

\bibitem{RoSt60}
G.~C. Rota and W.~G. Strang.
\newblock A note on the joint spectral radius.
\newblock {\em Indag. Math.}, 22:379--381, 1960.

\bibitem{Strum1999}
J.~F. Sturm.
\newblock Using \mbox{SeDuMi} 1.02, a \mbox{MATLAB} toolbox for optimization
  over symmetric cones.
\newblock {\em Optimization Methods and Software}, 11(12):625--653, 1999.
\newblock URL: http://sedumi.ie.lehigh.edu/.

\bibitem{BlTi3}
J.~N. Tsitsiklis and V.~Blondel.
\newblock The {L}yapunov exponent and joint spectral radius of pairs of
  matrices are hard- when not impossible- to compute and to approximate.
\newblock {\em Mathematics of Control, Signals, and Systems}, 10:31--40, 1997.

\end{thebibliography}


\end{document}